%% file: gl_pairing_final_II.tex
\documentclass[11pt]{amsart}
\usepackage{fullpage,xcolor,graphicx}
\usepackage[OT2,T1]{fontenc}
\usepackage{hyperref}
\usepackage{pifont}
\usepackage{amsthm}
\usepackage{amsfonts}
\usepackage{amssymb}
\usepackage[mathscr]{euscript}
\usepackage[all]{xy}
\usepackage{amsmath}
\usepackage{epsfig}
\usepackage{latexsym}
\usepackage{stackengine}
\usepackage{multirow}
\usepackage{hhline}
\usepackage{MnSymbol}

\newtheorem{theorem}{Theorem}[section]

\theoremstyle{definition}
\newtheorem{definition}[theorem]{Definition}
\newtheorem{example}[theorem]{Example}
\newtheorem*{acknowledgments*}{Acknowledgments}{}
\theoremstyle{remark}
\newtheorem{remark}[theorem]{Remark}

\newcommand{\lk}{\operatorname{\ell{\it k}}}

\makeatletter
\@namedef{subjclassname@2020}{\textup{2020} Mathematics Subject Classification}
\makeatother

\author[M. Chrisman]{Micah Chrisman}
\address{Department of Mathematics, The Ohio State University, Columbus, Ohio, USA}
\email{chrisman.76@osu.edu}

\dedicatory{In memory of Richard A. Litherland}

\begin{document}

\title[The G-L form and its many applications]{The Gordon-Litherland form and its many applications}

\begin{abstract} Gordon and Litherland's paper \emph{On the Signature of a link} introduced a bilinear form that simultaneously unifies both the quadratic forms of Trotter and Goeritz. This remarkable pairing of combinatorics and topology has had widespread application in low-dimensional topology.  In this expository note, we give a picture proof (via Kirby diagrams) of their main result and discuss the numerous ways their theorem has been put to good use.   
\end{abstract}


\maketitle

\section{Introduction}

This note is based on a talk given by the author in Autumn 2019 for the \emph{Invitations to Mathematics Seminar} at the Ohio State University. The purpose of this seminar is twofold. First, it introduces first-year graduate students to fundamental ideas which may otherwise be absent in their coursework. Secondly, the seminar aims to connect students with the various research communities in the Department of Mathematics. It is hard to find a topic more suited to these goals than the Gordon-Litherland pairing. Indeed, the main result of  \cite{gordon_litherland} succintly expresses the ideal to which much of modern knot theory aspires. Their theorem states that the signature of a knot $K$ is equal to the signature of its Goeritz matrix $G$ minus an error correction term $\mu(K)$:

\begin{theorem}[Gordon-Litherland] \label{thm_gl} $\sigma(K)=\text{sign}(G)-\mu(K)$.
\end{theorem}

The left-hand side of this equation is defined topologically whereas the right-hand side is defined combinatorially. In particular, both terms on the right can be quickly calculated from any diagram of the knot. But it is in the symbol ``$=$'' where the magic happens. It stands in place of a beautiful interplay of topology and diagrammatics that is nowadays a hallmark of the field of knot theory. The Gordon-Litherland form makes this connection explicit. By studying it, we gain insight into both classical topological tools (e.g. quadratic forms $\&$ branched covering spaces) and modern techniques (e.g. categorification $\&$ trisections). The Gordon-Litherland form thus serves as one of the best possible advertisements for the whole subject.

Other than the exposition, nothing in this note is new. The paper is geared towards early-career graduate students, so we present a proof of Theorem \ref{thm_gl} that avoids advanced topics such as obstruction theory and the $G$-signature theorem. Section \ref{sec_what} begins by defining each of the terms in Theorem \ref{thm_gl}. Our picture proof of Theorem \ref{thm_gl} is given in Section \ref{sec_pic_proof}. Applications are discussed in Section \ref{sec_app}. For generalizations of the Gordon-Litherland pairing, see Section \ref{sec_gen}. 

\section{What are $\sigma(K)$, $G$, and $\mu(K)$?} \label{sec_what}

\begin{figure}[htb]
\begin{tabular}{|ccc|} \hline
\begin{tabular}{c}  \includegraphics[scale=.3]{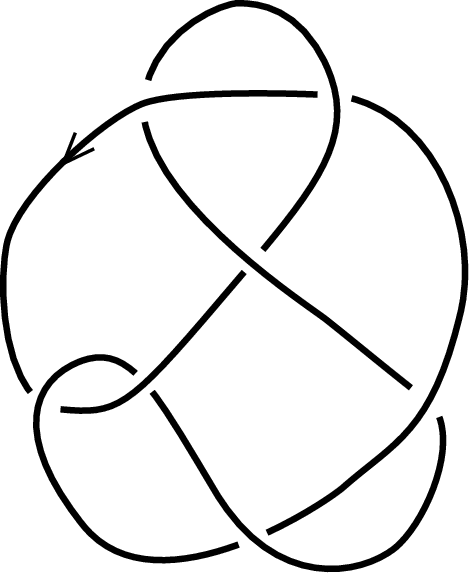} \\ \end{tabular} & 
\begin{tabular}{c}  \includegraphics[scale=.3]{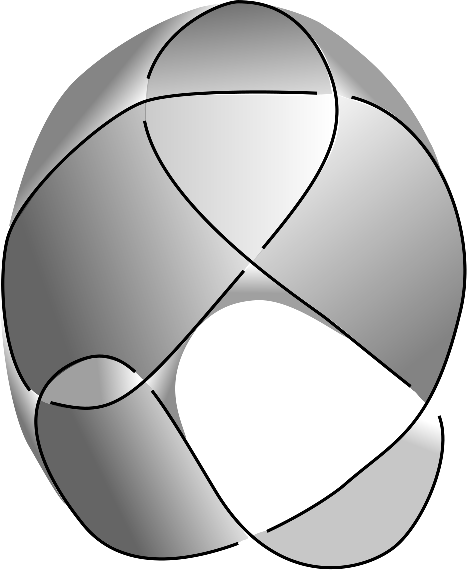} \\ \end{tabular} & 
\begin{tabular}{c} \\ \def\svgwidth{2.2in}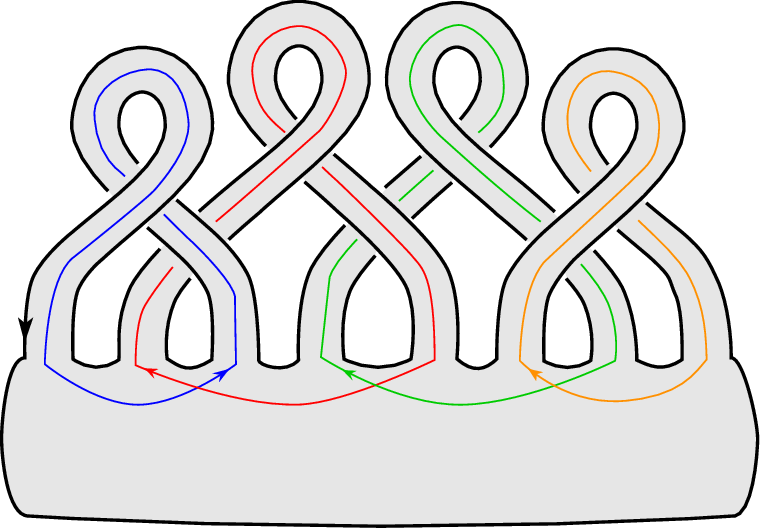 \\ \\ \end{tabular}  \\ \hline
\end{tabular}
\caption{The knot $7_6$, a Seifert surface for $7_6$, and the same Seifert surface isotoped into disc-band form showing the homology generators $\alpha_1,\alpha_2,\alpha_3,\alpha_4$.} \label{fig_7_6_seifert}
\end{figure}

For simplicity, we consider only the case that $K$ is a knot in $\mathbb{S}^3$. Let $F$ be a compact, connected, orientable surface in $\mathbb{S}^3$ such that $\partial F=K$. Then $F$ is called a \emph{Seifert surface} of $K$. A Seifert surface for the knot $7_6$ is given in Figure \ref{fig_7_6_seifert}. The Seifert pairing $\mathscr{S}_F:H_1(F;\mathbb{Z}) \times H_1(F;\mathbb{Z}) \to \mathbb{Z}$ is defined by $\mathscr{S}_F(\alpha,\beta)=\lk(\alpha,\beta^+)$, where $\beta^+$ denotes the positive push-off of $\beta$ into $\mathbb{S}^3 \smallsetminus F$. This simply means that since $F$ is oriented, it has a well-defined positive normal direction as determined by the right-hand rule. Hence, a closed curve $\beta$ on $F$ can be pushed in this direction to give a curve $\beta^+$ in $\mathbb{S}^3 \smallsetminus F$. Now, choose a basis $\{\alpha_1,\ldots,\alpha_{2g}\}$ of $H_1(F;\mathbb{Z})$ and let $A$ be the $2g \times 2g$ matrix whose $(i,j)$ entry is $\mathscr{S}_F(\alpha_i,\alpha_j)$. Then $A$ is called a \emph{Seifert matrix}. Denote by $A^{\intercal}$ the transpose of $A$. The \emph{signature} of $K$ is the signature of the symmetrized Seifert matrix:
\[
\sigma(K)=\text{sign}(A+A^{\intercal})
\] 
That is, $\sigma(K)$ is the number of positive eigenvalues of $A+A^{\intercal}$ minus the number of its negative eigenvalues. The signature of $K$ is independent of the choice of the Seifert surface $F$ and is an invariant of $K$.  

The Goeritz matrix of $K$ is defined as follows. Assume that $K$ is given by a knot diagram on $\mathbb{S}^2=\mathbb{R}^2 \cup \{\infty\} \subset \mathbb{S}^3$, which will also be denoted by $K$. The diagram separates $\mathbb{S}^2$ into regions and these regions have a checkerboard coloring, so that each edge borders a white region and a black region. We will assume that the checkerboard coloring is chosen so that $\infty$ lies in a white region. By the \emph{black surface} of $K$, we mean the spanning surface for $K$ made from the black regions of the checkerboard coloring together with half-twisted bands at the crossings. The \emph{white surface} is similarly defined. If $K$ is nontrivial, at least one of the black surface and the white surface is nonorientable (see Clark \cite{clark_78}, Theorem 2). Note that this is untrue in general for links. For example, both checkerboard surfaces of the usual Hopf link diagram are orientable.

\begin{figure}[htb]
\begin{tabular}{|cc|} \hline 
\begin{tabular}{c} \\ \tiny \def\svgwidth{1.2in}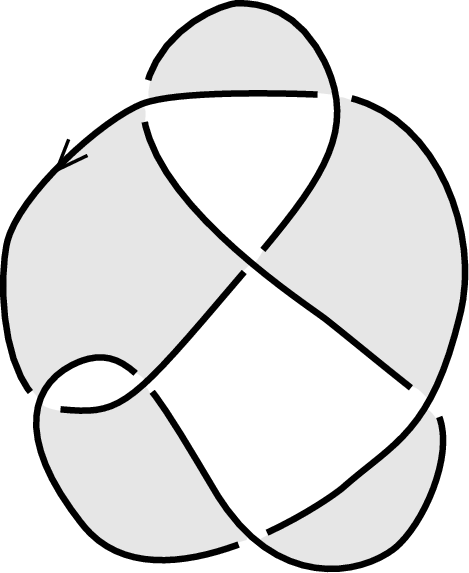 \\ \\  \end{tabular} & 
\begin{tabular}{c} \\ \def\svgwidth{2.2in}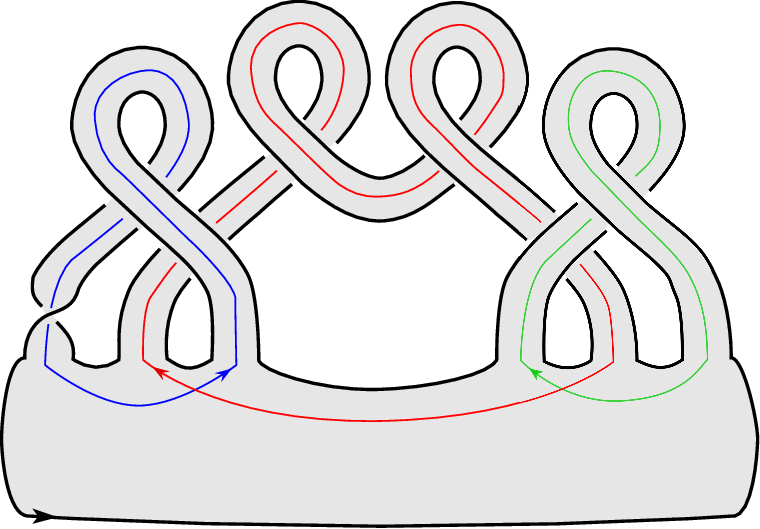 \\ \\ \end{tabular}  \\ \hline
\end{tabular}
\caption{A checkerboard coloring of $7_6$ and its black spanning surface isotoped into disc-band form.} \label{fig_7_6_spanning}
\end{figure}

\begin{figure}[htb]
\begin{tabular}{|cc|cc|} \hline
\begin{tabular}{c} \\ \def\svgwidth{.75in}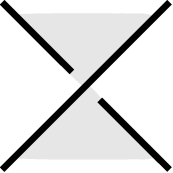 \\ $\eta(C)=1$ \\ \\ \end{tabular} & 
\begin{tabular}{c} \\ \def\svgwidth{.75in}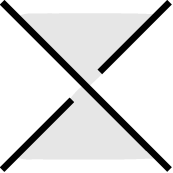 \\ $\eta(C)=-1$ \\ \\ \end{tabular} & 
\begin{tabular}{c} \\ \def\svgwidth{.75in}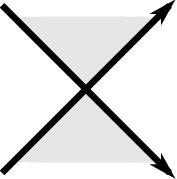 \\ Type I \\ \\ \end{tabular} &
\begin{tabular}{c} \\ \def\svgwidth{.75in}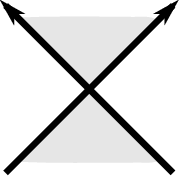 \\ Type II \\ \\ \end{tabular} \\ \hline 
\end{tabular}
\caption{Definitions of incidence number and type of a crossing.} \label{fig_incidence_type}
\end{figure}

Number the white regions of $K$ as $X_0,\ldots,X_N$. By convention, $X_0$ contains $\infty$. At each crossing $C$ of $K$, we have a an incidence number $\eta(C)$ as shown in Figure \ref{fig_incidence_type}.  For $0 \le i \ne j \le N$, define $S_{ij}$ to be the set of crossings of $K$ that are incident to both $X_i$ and $X_j$. Denote by $\delta_{ij}$ Kronecker's delta function and set:
\[
g_{ij}=\left\{\begin{array}{cl} -\sum_{C \in S_{ij}} \eta(C) & i \ne j \\ -\sum_{k=0}^N (1-\delta_{ik}) g_{ik} & i=j \end{array} \right.
\]
The Goeritz matrix is defined to be $G=(g_{ij})_{i,j=1}^N$. Unlike $\sigma(K)$, $\text{sign}(G)$ is not a knot invariant. That this is true will become clear in Section \ref{sec_step_4} ahead.

The necessary correction term $\mu(K)$ is defined as follows. A crossing $C$ of $K$ is said to be of \emph{type I} or \emph{type II} according to which of the pictures it matches on the right in in Figure \ref{fig_incidence_type}. Set $\mu(K)=\sum_{\text{type}(C)=II} \eta(C)$. 

\begin{example} \label{ex_7_6} Let $K=7_6$ and let $F$ be the Seifert surface shown in Figure \ref{fig_7_6_seifert}. On the right, $F$ has been deformed into a disc-band surface, so that its $0$-handles and $1$-handles are evident. The black surface for $K$ is shown in Figure \ref{fig_7_6_spanning}. Then we have: 
\[
A=\begin{pmatrix}
-1 & 0 & 0 & 0 \\
-1 & -1 & 0 & 0 \\
0 & 1 & 1 & -1 \\
0 & 0 & 0 & -1
\end{pmatrix}, \quad G=\begin{pmatrix}
3 & -1 & 0 \\
-1 & 4 & -1 \\
0 & -1 & 2
\end{pmatrix}, \quad \mu(K)=5.
\]
One can check that $A$ has 3 negative eigenvalues and $1$ positive eigenvalue whereas $G$ has $3$ positive eigenvalues. As expected, $\sigma(K)=\text{sign}(G)-\mu(K)$. 
\end{example}

\section{Illustrating the Proof} \label{sec_pic_proof}

\subsection{The Gordon-Litherland pairing} \label{sec_gl_pairing} The main idea of the proof is simple. First, we expand our domain of consideration to all \emph{spanning surfaces} $F$ of a knot $K$. That is, $F$ is any compact connected surface (orientable or not) such that $\partial F=K$. This includes both Seifert surfaces and checkerboard surfaces. Then we make everything orientable by passing to the twofold orientation cover $\widetilde{F}$ of $F$. The cover $\widetilde{F}$ will play an analogous role to that of the push-off in the Seifert form.

It is useful to construct $\widetilde{F}$ explicitly. Suppose that $F$ is given in disc-band form, that is as a union of $0$-handles (discs) and $1$-handles (bands). Examples are shown in Figures \ref{fig_7_6_seifert} and \ref{fig_7_6_spanning}. Both $0$-handles and $1$-handles are really just discs $\mathbb{B}^2$, and thus they lift as discs to the twofold cover. The cover $\widetilde{F}$ then inherits a handle decomposition from $F$; it is the union of the lifts of all $0$-handles and $1$-handles of $F$. Hence we may draw $\widetilde{F}$ by doubling all of the handles of $F$. For each $0$-handle $F^0$ of $F$, place duplicate $0$-handles $(F^0)',(F^0)''$ slightly above and below $F^0$, respectively, where ``above'' and ``below'' are relative to the direction of the knot projection. Similarly, each $1$-handle $F^1$ is duplicated to $(F^1)',(F^1)''$. Each copy of the band $F^1$ runs parallel to $F^1$ and its two ends are attached to the duplicates of some $0$-handles. Of course, a band that starts on a ``top'' duplicate $0$-handle may end on a ``bottom'' duplicate $0$-handle (and vice versa). This happens whenever $F^1$ has an odd number of half-twists. The level is preserved if $F^1$ has an even number of half-twists.

Given any closed path $\rho$ in $F$, let $\tau(\rho)$ denote the union of its two lifts in $\widetilde{F}$. Then either $\tau(\rho)$ consists of two closed paths or one closed path. That this is well-defined follows from the fact that the orientation cover is regular. Then $\tau$ descends to a well-defined map $\tau:H_1(F) \to H_1(\widetilde{F})$, which is called the \emph{transfer map}. 

\begin{definition}[Gordon-Litherland form] Let $K$ be a knot, $F$ any spanning surface for $K$, and $\tau:H_1(F) \to H_1(\widetilde{F})$ the transfer map. The \emph{Gordon-Litherland form} is the pairing $\mathcal{G}_F:H_1(F;\mathbb{Z}) \times H_1(F;\mathbb{Z}) \to \mathbb{Z}$ defined by:
\[
\mathcal{G}_F(\alpha,\beta)=\lk(\alpha,\tau(\beta)).
\] 
\end{definition}

The correction term $\mu(K)$ is realized geometrically by the \emph{Euler number}. This is defined as follows.  Denote by $\widehat{F}$ the surface $F$ pushed slightly into $\mathbb{B}^4$, so that $\partial F=K$ remains in $\mathbb{S}^3=\partial \mathbb{B}^4$. Let $\nu(\widehat{F})$ be the normal bundle of $\widehat{F} \subset \mathbb{B}^4$, $\widehat{F}'$ a non-trivial section of $\nu(\widehat{F})$, and $K'=\partial \widehat{F}'$. If $F$ is orientable, $\widehat{F}'$ is just a push-off of $\widehat{F}$ and it hence may be assumed that $\widehat{F} \cap \widehat{F}'=\varnothing$. But if $\widehat{F}$ is non-orientable, this cannot be assumed. A choice of local orientation of $\widehat{F}$ at each transversal intersection point $x \in \widehat{F} \cap \widehat{F}'$ induces a local orientation of $\widehat{F}'$ at $x$. The two pieces fit together to give a local orientation of $\mathbb{B}^4$. Hence, each $x \in \widehat{F} \cap \widehat{F}'$ has an incidence number $\pm 1$ corresponding to whether this orientation agrees or disagrees with the global orientation of $\mathbb{B}^4$. This is independent of the original choice of local orientation of $x \in \widehat{F}$. The sum of these incidence numbers is the \emph{Euler number} of $F$. Following the above conventions, it can be checked that:
\[
e(F)=-\lk(K,K').
\] 

The proof of Theorem \ref{thm_gl} from \cite{gordon_litherland} can now be outlined as follows. Below, $M_{\widehat{F}}$ denotes the twofold branched cover of $\mathbb{B}^4$ over $\widehat{F}$ and $\sigma(M_{\widehat{F}})$ denotes the signature of its intersection form. 
\begin{enumerate}
\item $\mathcal{G}_F$ recovers the symmetrized Seifert form when $F$ is a Seifert surface and the Goeritz form when $F$ is the black checkerboard surface.
\item The Euler number of the black surface is $-2 \mu(K)$.
\item The Gordon-Litherland form is isomorphic to the intersection form\footnote{For $F$ oriented, this was first proved by Kauffman and Taylor \cite{kauffman_taylor}.} on $M_{\widehat{F}}$.
\item Given any spanning surfaces $F_1,F_2$ of $K$,
\[
\sigma(M_{\widehat{F_1}})+\tfrac{1}{2}e(F_1)=\sigma(M_{\widehat{F_2}})+\tfrac{1}{2}e(F_2).
\]
\end{enumerate} 

With these established, the proof is then completed by letting $F_1$ be any orientable spanning surface of $K$, $F_2$ the black surface of $K$, and then applying the above formula:
\small
\[
\sigma(K) =\sigma(M_{\widehat{F_1}})+0=\sigma(M_{\widehat{F_1}})+\tfrac{1}{2}e(F_1)=\sigma(M_{\widehat{F_2}})+\tfrac{1}{2}e(F_2)=\text{sign}(G)+\tfrac{1}{2}(-2 \mu(K)).
\]
\normalsize It remains to prove each of the steps (1)-(4). This is done in the following subsections.

\subsection{Proof of Step (1)} Suppose that $F$ is a Seifert surface of $K$. Then $\widetilde{F}$ is just two copies of $F$. One copy is the positive push-off of $F$ from itself while the other copy is the negative push-off. Since the negative push-off is just the transpose of the positive push-off (see e.g. \cite{bz}, Lemma 8.6), it follows that $\mathcal{G}_F$ is the symmetrized Seifert form when $F$ is orientable.

\begin{figure}[htb]
\begin{tabular}{|c|cc|} \hline 
\begin{tabular}{c} \\  \def\svgwidth{1.75in}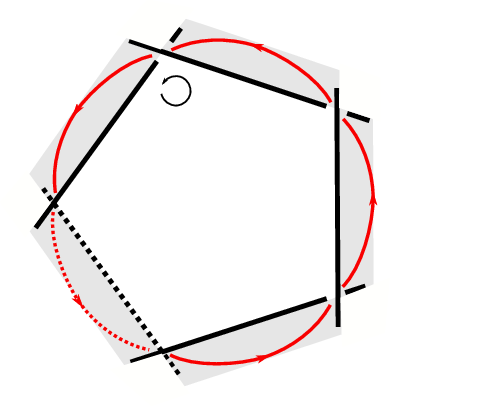 \\  \end{tabular} & 
\begin{tabular}{c} \\ \def\svgwidth{1.2in}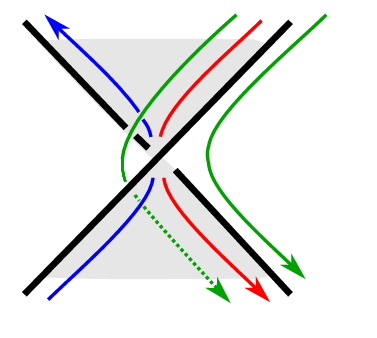 \\ $\eta(C)=1$\\ \\ \end{tabular}  & 
\begin{tabular}{c} \\ \def\svgwidth{1.2in}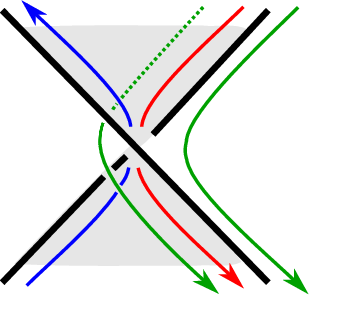 \\ $\eta(C)=-1$\\ \\ \end{tabular}  \\ \hline
\end{tabular}
\caption{Obtaining a basis for the black surface (left) and calculating the Gordon-Litherland form relative to that basis.} \label{fig_goeritz_mat_calc}
\end{figure}

Now suppose that $F$ is the black surface of $K$ and that the white regions are labeled $X_0,\ldots,X_N$. Orient each $X_i$ so that it matches the orientation of $\mathbb{S}^2$. Let $\alpha_i$ be a simple closed curve on $F$ that circumnavigates $X_i$ in the direction of its orientation (see Figure \ref{fig_goeritz_mat_calc}, left). The $1$-cycles $[\alpha_1],\ldots,[\alpha_N]$ form a basis for $H_1(F;\mathbb{\mathbb{Z}})$. This can be seen using an Euler characteristic argument as follows. Let $C$ be the number of crossings of the knot diagram and $B$ the number of black regions. A knot diagram is a planar graph with $C$ vertices, $2C$ edges, and $N+1+B$ faces. Then we have $C-2C+N+1+B=2$, so that $B=1+C-N$. The Euler characteristic of $F$ is $\chi(F)=1-\beta_1(F)$, where $\beta_1(F)$ is the first Betti number of $F$. A handle decomposition of $F$ is evident; each black region is a $0$-handle and and each half-twisted band at a crossing is a $1$-handle. Hence, $\chi(F)=B-C=1+C-N-C=1-N$ and it follows that $\beta_1(F)=N$. To see that $[\alpha_1],\ldots,[\alpha_N]$ are linearly independent, observe that $\sum_{i=0}^N [\alpha_i]$ bounds the union of the black discs. This implies that $\sum_{i=0}^N [\alpha_i]=0$. This is the only such linear relation, so discarding $[\alpha_0]$ then gives a basis for $H_1(F;\mathbb{Z})$.

It remains to identify the Goeritz matrix with $\mathcal{G}_F$. First suppose that $i \ne j$ and let $C$ be a crossing of $K$ and let $X_i,X_j$ the white regions incident to $C$. If $\eta(C)=1$, it follows from Figure \ref{fig_goeritz_mat_calc} that the contribution to $\lk(\alpha_i,\tau(\alpha_j))$ is $-1=-\eta(C)$. If $\eta(C)=-1$, the contribution is $1=-\eta(C)$. Hence $\mathcal{G}_F([\alpha_i],[\alpha_j])=g_{ij}$. If $i=j$, use $\sum_{i=0}^N [\alpha_i]=0$ to write $[\alpha_i]=-\sum_{j=0}^N (1-\delta_{ij})[\alpha_j]$. Then  $\mathcal{G}_F([\alpha_i],[\alpha_i])=\mathcal{G}_F([\alpha_i],-\sum_{j=0}^N (1-\delta_{ij})[\alpha_j])=g_{ii}$. This completes the proof of step (1).

\subsection{Proof of Step (2)} 
\begin{figure}[htb]
\begin{tabular}{|cccc|} \hline
\begin{tabular}{c} \\ \def\svgwidth{.75in}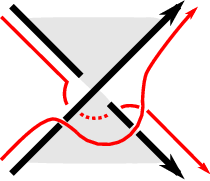 \\ $(1,I)$ \\ \\ \end{tabular} & 
\begin{tabular}{c} \\ \def\svgwidth{.75in}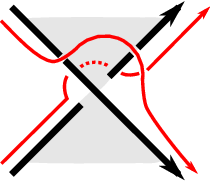 \\ $(-1,I)$ \\ \\ \end{tabular} & 
\begin{tabular}{c} \\ \def\svgwidth{.75in}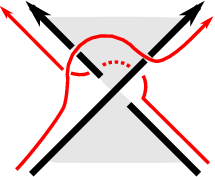 \\ $(1,II)$ \\ \\ \end{tabular} &
\begin{tabular}{c} \\ \def\svgwidth{.75in}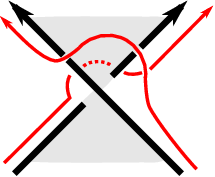 \\ $(-1,II)$ \\ \\ \end{tabular} \\ \hline 
\end{tabular}
\caption{The four possible pairs $(\eta(C),\text{type}(C))$ for a crossing $C$.} \label{fig_euler_calc}
\end{figure}

There are four possibilities for the pair $(\eta(C),\text{type}(C))$ where $C$ is a crossing of $K$. These are shown in Figure \ref{fig_euler_calc}, together with the curve $K'$. The contribution to $e(F)$ from type I crossings is zero and the contribution from type II crossings is $-2 
\cdot\eta(C)$. This completes the proof of step (2).

\subsection{Proof of Step (3)} If $M$ is a compact oriented $4$-manifold, the intersection form $Q_M:H^2(M,\partial M; \mathbb{Z}) \times H^2(M,\partial M; \mathbb{Z}) \to \mathbb{Z}$ is defined by $Q_M(\xi,\zeta)=\langle \xi \cup \zeta,[M]\rangle$, where $[M]$ is the fundamental class of $M$. The original argument for step (3) used a little algebraic topology and some properties of the intersection form. Here we will instead give a proof by picture. The idea is to draw the $4$-manifold $M$ and then calculate a matrix for $Q_M$ directly from this picture. The picture is an example of a \emph{Kirby diagram}. Kirby diagrams are a 4-dimensional analogue of the curve sketching techniques learned in first-semester calculus. They are a way of representing $4$-manifolds by drawing non-degenerate critical points of certain height functions.  First we will give an algorithm for drawing a Kirby diagram of the twofold branch cover of $\mathbb{B}^4$ with branching set $\widehat{F}$. Then we will briefly discuss why the algorithm works. This trick allows for a short proof of step (3), which is given at the end of this subsection. 

\begin{figure}[htb]
\begin{tabular}{|cccc|cccc|} \hline
& \begin{tabular}{c} \\ \def\svgwidth{.6in}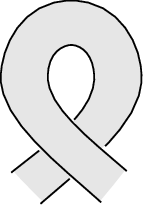 \\   \end{tabular} & \begin{tabular}{c} \\\Huge $\approx$ \normalsize \end{tabular} & 
\begin{tabular}{c} \\ \def\svgwidth{1in}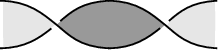 \\    \end{tabular} & 
\begin{tabular}{c} \\ \def\svgwidth{.6in}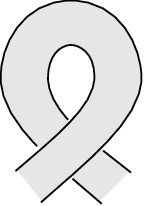 \\   \end{tabular} & \begin{tabular}{c} \\\Huge $\approx$ \normalsize \end{tabular} &
\begin{tabular}{c} \\ \def\svgwidth{1in}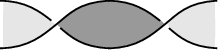 \\  \end{tabular} & \\ & \underline{a ``$+$'' curl:} & & \underline{two $+\tfrac{1}{2}$-twists:} & \underline{a ``$-$'' curl:} & & \underline{two $-\tfrac{1}{2}$-twists:} & \\ &  & & & & & & \\ \hline 
\end{tabular}
\caption{Calculating the sign of curls and half-twists.} \label{fig_curls_and_twists}
\end{figure}

\begin{figure}[htb]
\begin{tabular}{|ccc|} \hline 
\begin{tabular}{c} \\ \tiny \def\svgwidth{2.5in}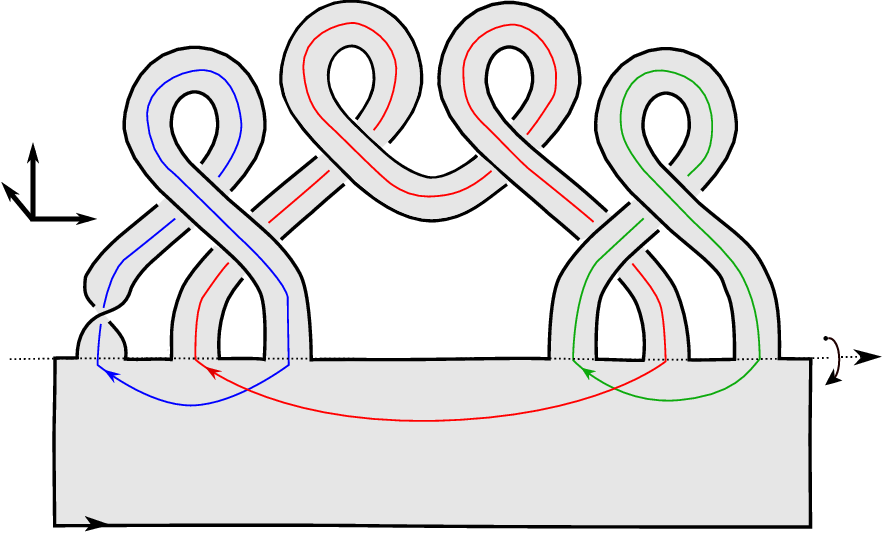 \\ \\  \end{tabular}  & \begin{tabular}{c} \\ \Huge $\longrightarrow$ \normalsize \\ \\ \end{tabular} & 
\begin{tabular}{c} \\ \def\svgwidth{2.2in}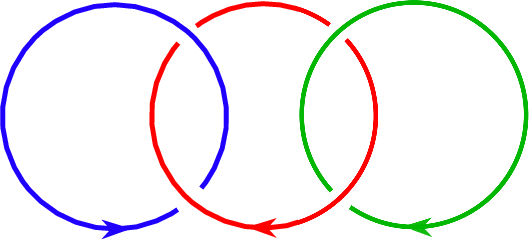 \\ \\ \end{tabular}  \\ \hline
\end{tabular}
\caption{Constructing the Kirby diagram for the double branch cover of $\mathbb{B}^4$ branched along the black surface $\widehat{F}$ of $7_6$.} \label{fig_7_6_double_cover}
\end{figure}

The algorithm we will use is due to Akbulut and Kirby \cite{akbulut_kirby}. Suppose that $F$ is given in disc-band form as a union of a single $0$-handle $F^0$ and $1$-handles $F^1_1,\ldots,F^1_n$. Hence,  $F=F^0 \cup F_1^1 \cup F_2^1 \cup \cdots \cup F_n^1$. Identify $\mathbb{S}^3$ with $\mathbb{R}^3 \cup \{\infty\}$. After an isotopy, it may be assumed that:
\begin{enumerate}
\item[($i$)] $F^0$ is a rectangle in the lower half-plane $\{(x,y,0)| y \le 0\}$ with one side on the $x$-axis, and
\item[($ii$)] each $F_i^1$ lies in the upper half-space $\{(x,y,z)| y \ge 0\}$ and meets the $xz$-plane only in the two intervals $f_i',f_i''$ along which it is attached to $F^0$. In particular, $f_i',f_i''$ lie on the $x$-axis. 
\end{enumerate}
See Figure \ref{fig_7_6_double_cover} for an example. Now, rotate the configuration of $1$-handles by an angle of $\pi$ about the $x$-axis. Then we have a second collection of $1$-handles $V_1^1 \cup \cdots \cup V_n^1$, with $V_i^1$ corresponding to $F_i^1$. The union of the cores of $F_i^1$ and $V_i^1$ is a knot $K_i$. If the band $F_i^1$ has $m_i$ half twists, then the annulus $V_i^1 \cup F_i^1$ has a total of $m_i$ full twists. The integer $m_i$ is called the \emph{framing} of $K_i$. Its sign is determined by the right-hand rule: a right-handed twist in a band is positive whereas a left-handed twist is negative (see Figure \ref{fig_curls_and_twists}). A knot $K_i$ with framing $m_i$ is denoted $K_i^{m_i}$. Define $L_F$ to be the framed link $K_1^{m_1} \cup \cdots \cup K_n^{m_n}$. In general, if $L=J_1^{m_1} \cup \cdots \cup J_n^{m_n}$ is a framed link, the \emph{linking matrix} of $L$ is the $n \times n$ matrix whose $(i,j)$ entry is $\lk(J_i,J_j)$  for $i \ne j$ and $m_i$ for $i=j$.    

\begin{example} Let $F$ be the nonorientable black surface for $7_6$, as shown in Figure \ref{fig_7_6_spanning}. A disc-band surface of $F$ is drawn in Figure \ref{fig_7_6_double_cover}, left. Using Figure \ref{fig_curls_and_twists}, we see that the first band has $3$ half-twists, the second band has $4$ half-twists, and the third band has $2$ half-twists. The Akbulut-Kirby algorithm returns the framed link $L_F$ on the right in Figure \ref{fig_7_6_double_cover}. It is quickly checked that the linking matrix for $L_F$ is exactly the Goeritz matrix $G$ from Example \ref{ex_7_6}. 
\end{example}

Why does this work? Briefly, a framed link $J_1^{m_1} \cup \cdots \cup J_n^{m_n}$ is a set of instructions for gluing $n$ $2$-handles $\mathbb{B}^2 \times \mathbb{B}^2$ to $\mathbb{B}^4$. For each $i$, it says to identify a small tubular neighborhood of $J_i$ with the solid torus $\mathbb{S}^1 \times \mathbb{B}^2=(\partial \mathbb{B}^2) \times \mathbb{B}^2$. The framing coefficient indicates how many times one should twist $\mathbb{S}^1 \times \mathbb{B}^2$ while gluing\footnote{To visualize this, take a bunch of spaghetti noodles, grab it with both hands, and start twisting.}. 

\begin{figure}[htb]
\begin{tabular}{|clc|} \hline 
\begin{tabular}{c} \\ \tiny \def\svgwidth{3in}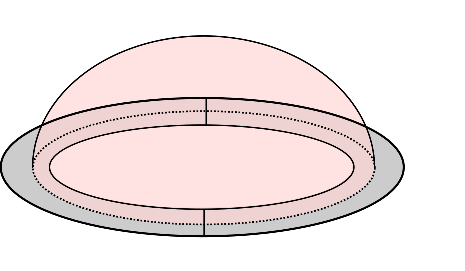 \\ \\  \end{tabular} & \hspace{-1cm} \begin{tabular}{l} \Huge $\rightarrow$\end{tabular} &
\begin{tabular}{c} \\ \includegraphics[scale=1]{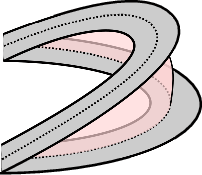} \\ \\ \end{tabular}  \\ \hline
\end{tabular}
\caption{Squishing the bubble gum so that the two 1-handles $F_i^1$, $V_i^1$ are identified.} \label{fig_bubblegum}
\end{figure}

The Kirby diagram for $M_{\widehat{F}}$ now arises as follows. First we construct the twofold branch cover of $\mathbb{B}^4$ with branching set $\widehat{F^0}$. With a little thought, one can convince oneself that cutting $\mathbb{B}^4$ along the disc $\widehat{F^0}$ gives $\mathbb{B}^4$ back again. However, we now have two copies of $F^0$, $F^0_+$ and $F^0_{-}$, with one copy on each side of the cut. Take two copies $B_1,B_2$ of this $\mathbb{B}^4$ and glue them together so that $F^0_+ \subset \partial B_1$ is identified with $F^0_-\subset B_2$ and $F^0_-\subset{B}_1$ is identified with $F^0_+ \subset B_2$. This is $M_{\widehat{F^0}}$, which is yet another $\mathbb{B}^4$. In $B_2$, we have the copies $V_1^1,\ldots,V_n^1$ of the $1$-handles $F_1^1,\ldots,F_n^1$. To get $M_{\widehat{F}}$, these must also be identified. This is done by thinking of a $2$-handle $\mathbb{B}^2 \times \mathbb{B}^2$ as a piece of $4$-dimensional bubble gum. For each $i$, the edge of the bubble gum is a solid torus $\mathbb{S}^1 \times \mathbb{B}^2$. The factor of $\mathbb{S}^1$ in this product, viewed as the unit circle in $\mathbb{C}$, is a union of two intervals: one $\mathbb{B}^1$ in the upper half-plane and another $\mathbb{B}^1$ in the lower half-plane. This splits $\mathbb{S}^1 \times \mathbb{B}^2$ into two halves, each of which is a $3$-ball $\mathbb{B}^1 \times \mathbb{B}^2$. Stick one of these halves along $F_i^1$, the other half along $V_i^1$, and then press $B_1$ and $B_2$ together so that the bubble gum gets squished. This is shown schematically in Figure \ref{fig_bubblegum}. The framing indicates the amount of twisting needed so that pressing $B_1$ and $B_2$ together aligns $F_i^1$ with $V_i^1$. The knot $\mathbb{S}^1 \times \{0\} \subset \mathbb{B}^2 \times \mathbb{B}^2$ is the knot $K_i$ in $L_F$. Performing this procedure for all $i$ gives our Kirby diagram.

Now we are almost ready to prove step (3). It remains only to compute the intersection form of $M_{\widehat{F}}$. This is easy to do from our Kirby diagram. We will summarize the standard argument for the reader's convenience (see e.g. \cite{gompf_stipsicz}, Proposition 4.5.11). By Poincar\'{e} duality, $H^2(M_{\widehat{F}},\partial M_{\widehat{F}}) \cong H_2(M_{\widehat{F}})$. A basis for $H_2(M_{\widehat{F}})$ can be read of the Kirby diagram. Let $\Sigma_i$ be a Seifert surface for $K_i$ and let $D_i$ by the disc $\mathbb{B}^2 \times \{0\} \subset \mathbb{B}^2 \times \mathbb{B}^2$ attached along $K_i$. Then it is not hard to believe that the homology classes $[D_1 \cup \Sigma_1],\ldots,[D_n \cup \Sigma_n]$ form a basis  $H_2(M_{\widehat{F}})$. In the duality isomorphism, the cup product is converted into an intersection of homology classes. This is fortuitous, as the intersections between $D_i \cup \Sigma_i$ and $D_j \cup \Sigma_j$ are visible in the Kirby diagram: they are the places where $K_j$ pierces $\Sigma_i$. Hence, the intersection form on $[D_i \cup \Sigma_i]$ and $[D_j \cup \Sigma_j]$ takes the value $\lk(K_i,K_j)$. It follows that the linking matrix is a matrix for $Q_{M_{\widehat{F}}}$. Thus, to prove step (3), it suffices to show that the linking matrix of $L_F$ is a matrix for $\mathcal{G}_F$. 

\begin{proof}[Proof of step (3)] The proof follows from several observations. Let $F$ be as above and let $\alpha_1,\ldots,\alpha_n$ be simple closed curves running along the $1$-handles, so that $[\alpha_1],\ldots,[\alpha_n]$ form a basis for $H_1(F)$. We will compute $\lk(\alpha_i,\tau(\alpha_j))$ by first adding up all crossing signs of $\alpha_i$ and $\tau(\alpha_j)$ and then dividing by $2$. If $\alpha_i,\alpha_j$ intersect in $F^0$, $\tau(\alpha_j)$ over-crosses and under-crosses $\alpha_i$ in the same direction. Hence, the contribution to $2 \cdot \lk(\alpha_i,\tau(\alpha_j))$ is $0$.  Next consider band crossings. The transfer map doubles the $\alpha_j$ band. Hence, if $i \ne j$ and $\alpha_i$ and $\alpha_j$ cross with sign $\varepsilon$, the contribution to $2 \cdot \lk(\alpha_i,\tau(\alpha_j))$ is $2 \varepsilon$. If $i=j$, both bands at the crossing are doubled and the contribution to $2 \cdot \lk(\alpha_i,\tau(\alpha_j))$ is $4 \varepsilon$. In particular, a $\pm$ curl contributes $\pm 4$. As a curl is two half-twists, a half-twist in a band contributes $\pm 2$. It follows that $\mathcal{G}_F(\alpha_i,\alpha_i)$ equals the number of half twists in the $i$-th band, which is the framing of $K_i$. For $i \ne j$, $\mathcal{G}_F(\alpha_i,\alpha_j)$ is the signed sum of the crossings of $F_i^1$ and $F_j^1$. By construction, each crossing of $K_i$ and $K_j$ in the upper half-plane occurs with the same sign as a crossing in the lower half-plane. Thus, $\lk(K_i,K_j)$ is also the signed sum of the band crossings and we conclude that $\mathcal{G}_F(\alpha_i,\alpha_j)=\lk(K_i,K_j)$. 
\end{proof}

\subsection{Proof of Step (4)} \label{sec_step_4}
\begin{figure}[htb]
\begin{tabular}{|c|ccc|} \hline  \rotatebox{90}{Tube Move \hspace{.8cm}} &
\def\svgwidth{2in}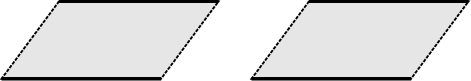 & \begin{tabular}{c} \\   \Huge $\leftrightarrow$ \normalsize \\ \\ \\ \end{tabular} & 
\tiny \def\svgwidth{2.5in}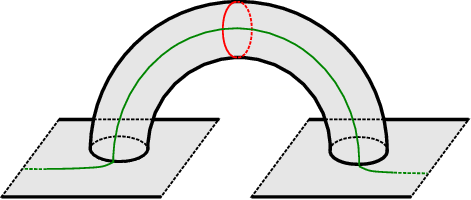 \normalsize \\ \hline & & & \\ \multirow[c]{2}*{\rotatebox{90}{Half-twist Moves\hspace{-.55in}}} & \multicolumn{3}{|c|}{\begin{tabular}{ccccc} 
\tiny \def\svgwidth{1.2in}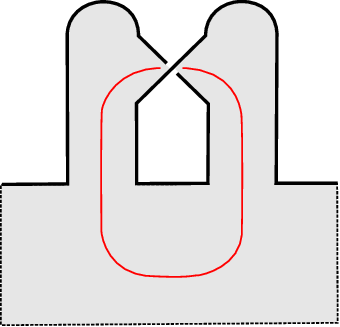 \normalsize & \begin{tabular}{c} \Huge $\leftrightarrow$ \normalsize \\ \\ \\ \end{tabular} &
 \def\svgwidth{1.2in}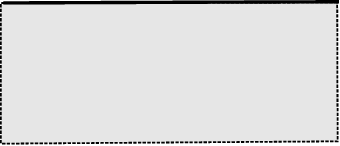 & \begin{tabular}{c} \Huge $\leftrightarrow$ \normalsize \\ \\ \\ \end{tabular} & 
\tiny \def\svgwidth{1.2in}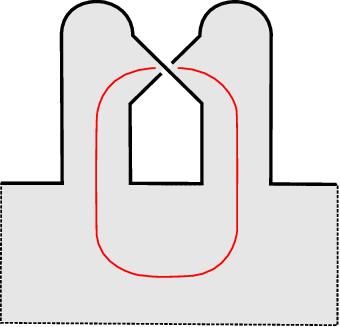 \normalsize \end{tabular}}  \\ \hline 
\end{tabular}
\caption{Isotopies, tube moves, and half-twist moves generate $S^*$-equivalence.} \label{fig_s_star}
\end{figure}

Two different arguments were given by Gordon and Litherland for the proof of step (4). One used the $G$-signature theorem (\cite{gordon_litherland}, Theorem 2) and the other relied upon obstruction theory (\cite{gordon_litherland}, Proposition 10 and Theorem 11). Fortunately, a theorem of Yasuhara allows us to stick with a more elementary approach. Two spanning surfaces of a link are said to be \emph{$S^*$-equivalent} if they may be obtained from one another by a finite sequence of ambient isotopies and the following operations:
\begin{enumerate}
\item Addition/deletion of a compressible tube (see Figure \ref{fig_s_star}, top ), and
\item Addition/deletion of a local half-twisted band (see Figure \ref{fig_s_star}, bottom).
\end{enumerate}
Here, a spanning surface need not be connected. Recall also that a tube $T$ is said to be \emph{compressible} if the meridian $m$ bounds an embedded disc $D$ in $\mathbb{S}^3$ that intersects $T$ only in $\partial D=m$. By analyzing how Reidemeister moves affect the checkerboard surfaces of a link, Yasuhara proved the following:

\begin{theorem}[Yasuhara \cite{yasuhara}] Any two checkerboard surfaces for a knot or link $K$ are $S^*$-equivalent and any spanning surface of $K$ is $S^*$-equivalent to a checkerboard surface for $K$.
\end{theorem}

\begin{figure}[htb]
\begin{tabular}{|ccc|} \hline  
\begin{tabular}{c} \tiny \def\svgwidth{1.9in}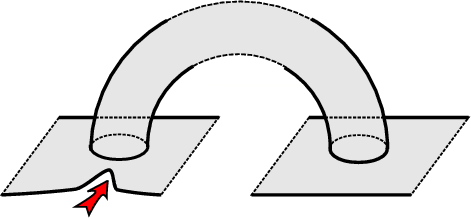 \normalsize \\ $t=0$ \\ \end{tabular} & \begin{tabular}{c} 
 \def\svgwidth{1.9in}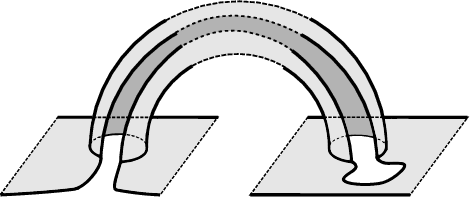 \\ $t=\tfrac{1}{2}$ \\ \\ \end{tabular} & \begin{tabular}{c}  \\ \\
\tiny \def\svgwidth{1.9in}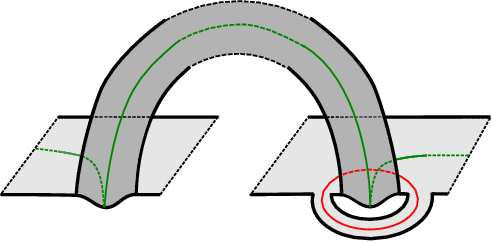 \normalsize \\ $t=1$\\ \end{tabular} \\ \hline \end{tabular}
\caption{Isotopy of a compressible tube to disc-band form.} \label{fig_iso}
\end{figure}

Thanks to Yasuhara's theorem, we can persevere with our proof by Kirby calculus. By step (3), it is only necessary determine the effect of $S^*$-equivalence on the linking matrix of the framed link $L_F$ representing $M_{\widehat{F}}$. Let $G$ be a matrix for $\mathcal{G}_F$ with respect to some basis $\{\alpha_1,\ldots,\alpha_n\}$ for $H_1(F)$.  An isotopy of $F$ clearly does not change $M_{\widehat{F}}$ and hence the only affect on $G$ is a unimodular change of basis. The Euler number of $F$ is likewise unaffected by ambient isotopy. Now, a half-twist move adds or removes a disjoint unknotted component to the Kirby diagram of $M_{\widehat{F}}$ with framing $\pm 1$: $L_F \leftrightarrow L_F \sqcup \bigcirc^{\pm 1}$ (see Figure \ref{fig_s_star}, bottom). The effect on the linking matrix is:
\[
G \leftrightarrow \left(\begin{array}{c|c} G & \textbf{0} \\ \hline \textbf{0} & \pm 1 \end{array}\right),
\]
From Figure \ref{fig_euler_calc}, it follows that a right-handed twist in a band contributes $-2$ to the Euler number and a left-handed twist contributes $+2$. Hence, $\sigma(M_{\widehat{F}})+\frac{1}{2}e(F)$ is unchanged by addition or removal of a local half-twisted band. Lastly, consider the effect of adding a compressible tube to $F$. Since $K$ is a knot, we may assume $F$ is connected. Then a tube move adds two new generators $\alpha_{n+1},\alpha_{n+2}$ to first homology (see Figure \ref{fig_s_star}, top left). After an isotopy (see Figure \ref{fig_iso}) and performing the Akbulut-Kirby algorithm, two components $K_{n+1}$, $K_{n+2}$ are added to $L_F$. Observe that $\lk(K_i,K_{n+2})=0$ for $1 \le i \le n$ and $\lk(K_{n+2},K_{n+1})=\pm 1$. The effect on the linking matrix is:
\[
G \to \left(\begin{array}{c|c|c} G & b & \textbf{0} \\ \hline b^{\intercal} & a & \pm 1 \\ \hline \textbf{0} & \pm 1 & 0 \end{array} \right),
\]
where $b$ is an $n \times 1$ vector and $a \in \mathbb{Z}$. Using simultaneous row and column operations, the above matrix is unimodular congruent to:
\[
\left(\begin{array}{c|c|c} G & \textbf{0} & \textbf{0} \\ \hline \textbf{0} & c & \pm 1 \\ \hline \textbf{0} & \pm 1 & 0 \end{array} \right),
\]
for some integer $c$. The lower $2\times 2$ block has a positive eigenvalue $\frac{1}{2} \left(c+\sqrt{c^2+4}\right)$ and a negative eigenvalue $\frac{1}{2} \left(c-\sqrt{c^2+4}\right)$. It follows that the signature is unchanged by a tube move. Since $e(F)=-\lk(K,K')$ and tube moves do not affect the boundary of $F$, the Euler number is also unchanged by a tube move. This completes the proof of step (4) and Theorem \ref{thm_gl}. \hfill $\square$

\begin{remark} To prove Theorem \ref{thm_gl}, it is only necessary to prove step (4) for spanning surfaces of a knot $K$. However, step (4) holds for any smoothly and properly embedded connected surface $S \subset \mathbb{B}^4$ with $\partial S=K$. This is the part of \cite{gordon_litherland} that uses the $G$-signature theorem.
\end{remark}

\section{Applications} \label{sec_app}

Since its appearance in 1978, the Gordon-Litherland form has found many uses in knot theory and beyond (see e.g. \cite{alexander_machon,nicholson}). Here we will focus on how the Gordon-Litherland form can be used to establish practical bounds on geometric properties of knots. Five instances of this theme are explored below.

\subsection{Gordian distances} 
\begin{figure}[htb]
\begin{tabular}{|ccc|} \hline & & \\
\begin{tabular}{c}
\def\svgwidth{.75in}
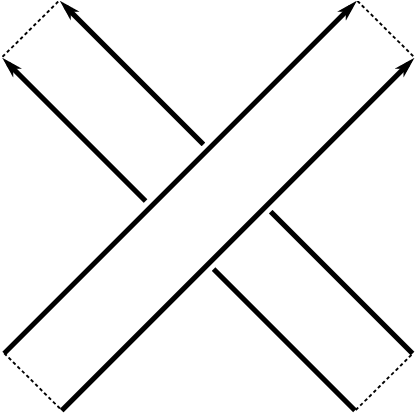 \end{tabular} & \begin{tabular}{c} \Huge $\leftrightarrow$ \end{tabular} & \begin{tabular}{c} \def\svgwidth{.75in}
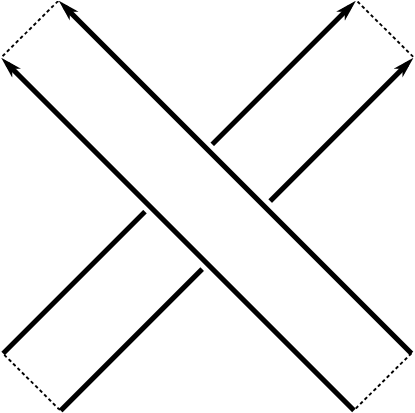 
\end{tabular} \\ & & \\ \hline
\end{tabular}
\caption{The $\sharp$-move.}\label{fig_sharp}
\end{figure}

\begin{figure}[htb]
\begin{tabular}{|ccc|} \hline & & \\
\begin{tabular}{c}
\includegraphics[scale=.4]{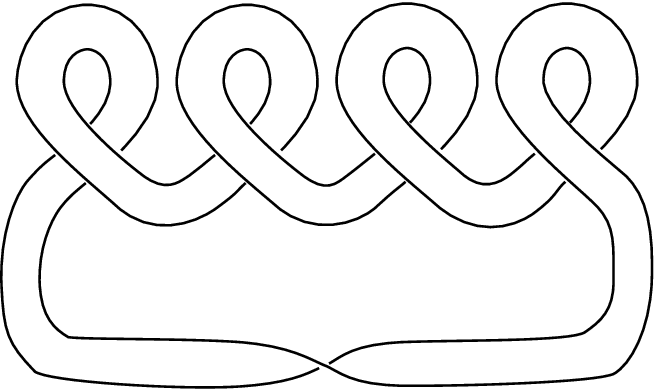} \end{tabular} & \begin{tabular}{c} \Huge $\rightarrow$ \end{tabular} & \begin{tabular}{c} 
\includegraphics[scale=.4]{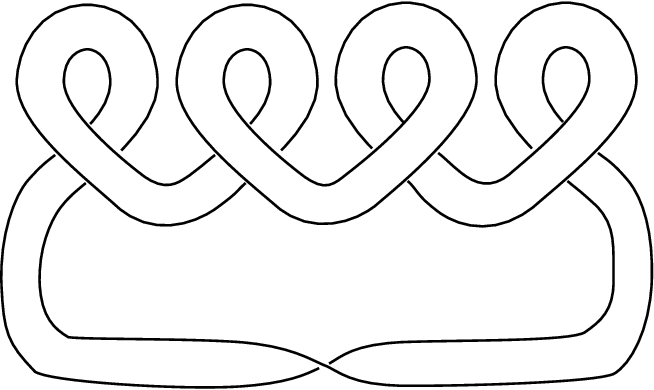} 
\end{tabular} \\ & & \\ \hline
\end{tabular}
\caption{Two $\sharp$-moves convert $T(9,2)$ (left) to the unknot (right).}\label{fig_t_9_2}
\end{figure}

If $K_1,K_2$ are oriented knots, the \emph{Gordian distance} is the minimum number of crossing changes needed to take some diagram of $K_1$ to some diagram of $K_2$. This is denoted by $d_G(K_1,K_2)$.  If $\bigcirc$ denotes the unknot, then $d_G(K,\bigcirc)=u(K)$ is the unknotting number of $K$. This defines a metric on the set of oriented knot types. Since $|\sigma(K_1)-\sigma(K_2)| \le 2$ when $K_2$ is obtained from $K_1$ by a crossing change, it follows that:
\[
d_G(K_1,K_2) \ge \frac{1}{2}|\sigma(K_1)-\sigma(K_2)|.
\]
Another uknotting operation is the $\sharp$-move (Figure \ref{fig_sharp}). The $\sharp$-\emph{Gordian distance} $d^{\sharp}_G(K_1,K_2)$ is the minimum number of $\sharp$-moves needed to convert a diagram of $K_1$ to $K_2$. The $\sharp$-\emph{unknotting number} of $K$ is $u^{\sharp}(K)=d_{\sharp}(K,\bigcirc)$. In \cite{murakami}, Murakami used the Gordon-Litherland pairing to obtain the following lower bound on $d^{\sharp}_G$:
\[
d^{\sharp}_G(K_1,K_2) \ge \frac{1}{6}|\sigma(K_1)-\sigma(K_2)|.
\]
Consider, for example the $T(m,2)$ torus knot for $m=2k+1$ an odd positive integer. Clearly, an upper bound on $u(K)$ is $k$, which is is obtained by switching any $k$ of its crossings. Using the Gordon-Litherland pairing, we easily calculate that $\sigma(T(m,2))=-m+1=2k$. Hence we have $k \le u(K) \le k$ and $u(K)=k$. On the other hand, the $\sharp$-unknotting number can be less than the unknotting number. Indeed, $u(T(9,2))=4$ but $u^{\sharp}(T(9,2))=2$. Murakami's lower bound gives $u^{\sharp}(T(9,2)) > 4/3$, so that $u^{\sharp}(T(9,2)) \ge 2$. The lower bound of $2$ is realized as follows. Covert each adjacent pair of crossings in $T(9,2)$ into full twists (see Figure \ref{fig_t_9_2}). This gives $4$ full twists in the same direction with one crossing left over. Since opposite curls cancel, performing $\sharp$-moves on any two of the curls gives an unknot. Hence, $2 \le u^{\sharp}(T(9,2))\le 2 $.  Further bounds on $d_G^{\sharp}$ in terms of the Alexander module were found by Murakami-Sakai \cite{murakami_sakai}. The $\sharp$-unknotting numbers of many torus knots were computed by Kanenobu \cite{kanenobu}. 

\subsection{Crosscap numbers} The \emph{$3$-genus} $g(K)$ of a knot $K$ is the smallest genus of any of its Seifert surfaces. The \emph{crosscap number} $C(K)$ of a knot $K$ in $\mathbb{S}^3$ is the smallest first Betti number among all nonorientable surfaces that span $K$. By definition, the crosscap number of the unknot is $0$. Crosscap numbers were first studied by Clark \cite{clark_78}, in the same year that \cite{gordon_litherland} was published. A knot with crosscap number $1$ bounds a M\"{o}bius band in $\mathbb{S}^3$. Clark proved that the only knots bounding M\"{o}bius bands are $(2,n)$ cable knots. In particular, every hyperbolic knot must have crosscap number at least $2$. As the trefoil is the $(3,2)$ torus knot, it has crosscap number $1$. For the figure eight knot $4_1$, which is hyperbolic, it is easy to construct a $\beta_1=2$ spanning surface from the standard projection. It follows that $C(4_1)=2$. Clark also showed that $C(K) \le 2 \cdot g(K)+1$ for any knot $K$ and asked if this upper bound was ever achieved.

Clark's question was finally answered in 1995 by Murakami and Yasuhara \cite{murakami_yasuhara} using the Gordon-Litherland pairing. To do this, they first proved the following necessary condition for a knot $K$ to bound a $\beta_1=2$ nonorientable spanning surface $F$.

\begin{theorem}[Murakami-Yasuhara \cite{murakami_yasuhara}] \label{murakami_yasuhara} There is a pair of curves generating $H_1(F;\mathbb{Z})$ such that the Goeritz matrix of $F$ takes the form:
\[
G_F(K)=\begin{bmatrix}
l & m \\ m & n
\end{bmatrix}
\] 
and satisfies $\sigma(K)=\text{sign}(G_F(K))-(l+2m+n)$, where $l,n$ are odd and $m$ is even.
\end{theorem} 

Consider, as in \cite{murakami_yasuhara}, the knot $7_4$. This knot has $3$-genus $1$ and hence $C(7_4) \le 3$. Since $7_4$ is hyperbolic, we have $C(7_4) \ge 2$. It therefore suffices to show that $C(K) \ne 2$. Now, suppose that $C(K)=2$. Since $7_4$ has determinant 15, one can write out integral binary quadratic forms with determinant 15 (see e.g. Conway-Sloane \cite{conway_sloane}). Using a case analysis, it is then possible to show that no such binary quadratic form satisfies Theorem \ref{murakami_yasuhara}. Hence, it follows that $C(7_4)=3=2 \cdot g(7_4)+1$. 

Crosscap numbers are in general difficult to compute. Yet progress is currently being made. For example, Kalfagianni-Lee \cite{kalfagianni_lee} determined upper and lower bounds on the crosscap numbers for several families of knots in terms of the coefficients of the Jones polynomial. 

\subsection{Non-orientable slice genus} The \emph{(smooth) 4-genus} of a knot $K$ is the smallest genus among all surfaces properly embedded in $\mathbb{B}^4$ such that $\partial F=K$. A well-known lower bound on the $4$-genus is $|\sigma(K)|/2$. In analogy with crosscap number, the \emph{(smooth) non-orientable slice genus} $\gamma^*(K)$ of a knot $K$ in $\mathbb{S}^3$ is defined as the smallest first Betti number of all compact connected non-orientable surfaces $F$ smoothly embedded in $\mathbb{B}^4$ such that $\partial F=K$. If $K$ is slice, then $\gamma^*(K)$ is defined to be $0$. In 1975, Viro \cite{viro} showed that $\gamma^*(4_1)=2$, and hence $4_1$ does not bound a M\"{o}bius band in $\mathbb{B}^4$. 

Another proof that $\gamma^*(4_1)=2$ follows from a result of Yasuhara, which is proved using the Gordon-Litherland pairing. Recall that the quadratic form $q:H_1(F;\mathbb{Z}/2\mathbb{Z}) \to \mathbb{Z}/2\mathbb{Z}$ of a knot $K$ is defined by $q(x)=x^{\intercal} A x$, where $F$ is a Seifert surface of $K$ and $A$ is a Seifert matrix for $F$. The \emph{Arf invariant} of $K$ is the Arf invariant of this quadratic form. Hence, $\text{Arf}(K)=0$ if $q$ maps a majority of elements of $H_1(F;\mathbb{Z}/2\mathbb{Z})$ to $0$ and otherwise $\text{Arf}(K)=1$. The Arf invariant is also the $\!\!\pmod 2$ coefficient of $z^2$ in the Conway polynomial of $K$. Note that there is a well-defined injection $\mathbb{Z}/2\mathbb{Z} \to \mathbb{Z}/8 \mathbb{Z}$ given by $x \to 4 \cdot x$.

\begin{theorem}[Yasuhara \cite{yasuhara_96}] If a knot $K$ in $\mathbb{S}^3$ bounds a M\"{o}bius band in $\mathbb{B}^4$, then
\[
\sigma(K)+4 \cdot \text{Arf}(K)  \equiv 0 \text{ or } \pm 2 \pmod 8.
\]
\end{theorem} 

Here we use the statement of Yasuhara's theorem given by Gilmer-Livingston \cite{gilmer_livingston}, to which we refer the reader for a short proof. Since $\sigma(4_1)=0$ and $\text{Arf}(4_1)=1$, it follows from Yasuhara's theorem that $\gamma^*(4_1)\ge 2$. Since $4_1$ bounds a crosscap number $2$ surface in $\mathbb{S}^3$, this surface can pushed into $\mathbb{B}^4$ to realize the lower bound. Hence $\gamma^*(4_1)=2$. By a similar argument, it follows that $\gamma^*(3_1 \sharp 3_1)=2$. 

The Gordon-Litherland pairing again plays a fundamental role in obstructing knots from bounding punctured Klein bottles. This is due to the following theorem of Gilmer and Livingston.

\begin{theorem}[Gilmer-Livingston \cite{gilmer_livingston}] If $K$ is a knot, $F \subset \mathbb{B}^4$ is a punctured Klein bottle with $\partial F=K$, and the $2$-fold branched cover $M_F$ over $F$ has a positive definite intersection form, then $\sigma(K)+4 \cdot \text{Arf}(K) \equiv 0,2, \text{ or } 4 \pmod 8$. On the other hand, if $M_F$ has a negative definite intersection form, $\sigma(K)+4 \cdot \text{Arf}(K) \equiv 0,4, \text{ or } 6 \pmod 8$. 
\end{theorem}

This result, together with some tools from Heegard-Floer homology, was used to prove that there are knots $K$ satisfying $\gamma^*(K) \ge 3$.  More recently, Batson \cite{batson} used Heegard-Floer homology to show that the the nonorientable slice genus of a knot can be arbitrarily large. In  \cite{jabuka_vancott}, Jabuka and Van Cott showed that the difference between the crosscap number and the nonorientable slice genus can be arbitrarily large for torus knots. The nonorientable $4$-genera of many double twist knots were determined by Hoste, Shanahan, and Van Cott \cite{hsvc}. 

\subsection{Alternating links $\&$ categorification} For each link $L$, Khovanov \cite{khovanov} introduced a bigraded chain complex whose cohomology groups $\text{Kh}^{i,j}(L)$ give a link invariant. The graded Euler characteristic is the Jones polynomial and hence Khovanov homology is said to \emph{categorify} the Jones polynomial. Knot Floer homology, independently discovered by Ozsv\'{a}th-Szab\'{o} \cite{os} and Rasmussen \cite{rasmussen_thesis}, also associates to each knot $K$ in $\mathbb{S}^3$ a bigraded chain complex whose homology $\widehat{HFK}(K)$ is a knot invariant. In this case, the graded Euler characteristic of $\widehat{HFK}(K)$ is the Alexander-Conway polynomial. 

An interesting application of Theorem \ref{thm_gl} occurs in the calculation of Khovanov homology and knot Floer homology for alternating links. In \cite{lee}, Lee proved that a non-split alternating link has its Khovanov homology with rational coefficients completely determined its Jones polynomial and signature. This was later shown to also be true in the case of integer coefficients by Shumakovitch \cite{shumakovitch_1,shumakovitch_2}. In the Heegard-Floer setting, Ozsv\'{a}th and Szab\'{o} proved that the knot Floer homology of an alternating knot is completely determined by its signature and its Alexander-Conway polynomial. The proof in both situations rely on the special form Theorem \ref{thm_gl} takes on alternating knots:
\[
\sigma(K)=\#(\text{black regions})-\#(\text{positive crossings})-1.
\]

Knot homology theories, together with Theorem \ref{thm_gl}, can be used to estimate how far a knot is away from being alternating. This is done using the \emph{Turaev genus} \cite{turaev_87}. Given a diagram $D$ of a knot $K$, there is an embedded oriented surface $\Sigma_D \subset \mathbb{S}^3$ such that $K$ projects to an alternating knot diagram $D'$ on $\Sigma_D$. The minimum genus of $\Sigma_D$ over all diagrams $D$ of the knot type is called the \emph{Turaev genus} $g_T(K)$ of $K$. Hence, $g_T(K)=0$ if and only if $K$ is alternating. The higher the Turaev genus, the further a knot is from being alternating. Lower bounds on the Turaev genus can be found in terms of the Rasmussen invariant $s(K)$ and the $\tau$-invariant $\tau(K)$. The Rasmussen invariant is extracted from Khovanov homology whereas the $\tau$-invariant is calculated from Heegard-Floer homology. In \cite{dasbach_lowrance}, Dasbach and Lowrance prove the following lower bounds on the Turaev genus:
\[
\left|\tau(K)+ \frac{\sigma(K)}{2} \right| \le g_T(K), \quad \left|\frac{s(K)+\sigma(K)}{2} \right| \le g_T(K), \text{ and } \left|\tau(K)- \frac{s(K)}{2} \right| \le g_T(K).
\]
When $K$ is alternating, this recovers the well-known relationship between $s$, $\tau$, and $\sigma$: $s(K)=2 \cdot \tau(K)=-\sigma(K)$. Using the inequalities above, Dasbach and Lowrance were able to determine the Turaev genus exactly for several families of $3$-braid torus knots.


\subsection{The Massey-Whitney Theorem} As a final application, it is worthwhile to mention the Massey-Whitney Theorem. In 1940, Whitney conjectured that the Euler number of a closed nonorientable surface $S$ smoothly embedded in $\mathbb{S}^4$ is bounded in absolute value by $4-2\chi$, where $\chi$ is the Euler characteristic of $S$. This was first proved by Massey \cite{massey} in 1969 using the $G$-signature theorem:

\begin{theorem}[Massey-Whitney] \label{thm_m_w} If $S$ is a closed, connected, nonorientable surface smoothly embedded in $\mathbb{S}^4$ with Euler characteristic $\chi$, then its Euler number is one of the following integers:
\[
2 \chi-4,\quad 2 \chi,\quad 2 \chi+4,\ldots,-2\chi-4,-2\chi, 4-2\chi.
\]
Moreover, each of these integers is realizable by some embedding of $S$ into $\mathbb{S}^4$.
\end{theorem}

A new proof of the Massey-Whitney Theorem was obtained by Joseph, Meier, Miller, and Zupan \cite{JMMZ} using the theory of trisections. Recall that every closed $3$-manifold $M$ has a \emph{Heegard splitting}. That is, $M$ can be decomposed into two $3$-dimensional $1$-handlebodies which are identified along their boundary surfaces. Gay and Kirby generalized this idea to $4$-manifolds \cite{gay_kirby}. They showed that every closed $4$-manifold $M$ can be written as $M=M_1 \cup M_2 \cup M_3$, where each $M_i$ is a $4$-dimensional $1$-handlebody, $M_i \cap M_j$ is a $3$-dimensional handlebody for $i \ne j$, and $M_1 \cap M_2 \cap M_3$ is a closed surface. Meier and Zupan \cite{meier_zupan_2} proved that every smoothly embedded surface $S$ in a trisected $4$-manifold admits a \emph{bridge position}, which means that $S \cap M_i$ is a collection of unknotted discs for each $i$ and $S \cap (M_i \cap M_j)$ is a trivial tangle for $i \ne j$. Bridge trisections thus bring knotted surfaces in $\mathbb{S}^4$ down to a level where which they can be studied using techniques from classical knot theory. 

For example, in the standard genus zero trisection of $\mathbb{S}^4$, every knotted surface has a tri-plane diagram $\mathbb{D}=(\mathbb{D}_1,\mathbb{D}_2,\mathbb{D}_3)$, where $\mathbb{D}_i$ is a trivial tangle diagram and $\mathbb{D}_i \cup \overline{\mathbb{D}}_{i+1}$ is a planar unlink. Then the Euler number can be recovered as the sum of the writhes of the diagram $\mathbb{D}_1 \cup \overline{\mathbb{D}}_{2}$, $\mathbb{D}_2 \cup \overline{\mathbb{D}}_{3}$, $\mathbb{D}_3 \cup \overline{\mathbb{D}}_{1}$. The proof of Theorem \ref{thm_m_w} follows from this fact and several applications of the Gordon-Litherland form.

\section{Generalizations} \label{sec_gen}

Owing to its success in classical knot theory, it is natural to hope that the Gordon-Litherland form might extend to knots in other $3$-manifolds. However, an immediate difficulty arises. If $M$ is a $3$-manifold such that $H_1(M;\mathbb{Z}) \ne 0$, then not all knots will bound Seifert surfaces and if $H_1(M;\mathbb{Z}/2\mathbb{Z}) \ne 0$, then not all knots in $M$ will bound spanning surfaces. Generalizations of the Gordon-Litherland form therefore either make restrictions on knots in $M$ or on $M$ itself so that these homological obstructions are removed.    

\subsection{Knots in $\mathbb{Z}/2\mathbb{Z}$-homology spheres} \label{sec_z_2z} The Gordon-Litherland pairing was generalized to knots in  $\mathbb{Z}/2\mathbb{Z}$-homology spheres by Greene \cite{greene}. To describe this generalization, it is first necessary to define the correct notion of linking number. Let $Y$ be a $\mathbb{Z}/2\mathbb{Z}$-homology sphere and let $K_1,K_2$ be knots in $Y$. By the universal coefficient theorem, $Y$ is also a rational homology $3$-sphere. Hence, there is an integer $q$ such that $q [K_1]=0 \in H_1(Y;\mathbb{Z})$ and a surface $S$ in the exterior of $K_1$ that wraps $q$ times around a longitude of $K_1$. The \emph{rational linking number} is defined to be $\lk_{Y,\mathbb{Q}}(K_1,K_2)=S \bullet K_2/q$, where $\bullet$ denotes the algebraic intersection number. Like the classical linking number, the rational linking number is well-defined and symmetric. 

Now if $F$ is a spanning surface for a knot $K$ in $Y$, a \emph{rational Gordon-Litherland form} $\mathcal{G}_F:H_1(F;\mathbb{Z}) \times H_1(F;\mathbb{Z}) \to \mathbb{Q}$ is defined by:
\[
\mathcal{G}_F(\alpha,\beta)=\lk_{Y,\mathbb{Q}}(\alpha,\tau(\beta)),
\]
where $\tau$ again denotes the transfer map. As in the classical case, the rational Gordon-Litherland form $\mathcal{G}_F$ is symmetric. Greene shows that a version of Theorem \ref{thm_gl} also holds for links in $\mathbb{Z}/2\mathbb{Z}$-homology $3$-spheres. If $L=K_1 \cup K_2 \cup \cdots \cup K_n$ is an unoriented link, the quantity $\text{sign}(\mathcal{G}_F)+\tfrac{1}{2}e(F)$ is the same for all spanning surfaces $F$ of $L$.  Here $e(F)$ is the Euler number of $F$, which is defined analogously to the classical case given in Section \ref{sec_gl_pairing}. When $Y=S^3$, $\text{sign}(\mathcal{G}_F)+\tfrac{1}{2}e(F)$ is the \emph{Murasugi signature}, which is the average of the signatures of the oriented links with underlying unoriented link $L$.  If $L$ is oriented, $\text{sign}(\mathcal{G}_F)+\tfrac{1}{2}e(F,L)$ is also independent of the choice of $F$, where $e(F,L)=e(F)-\sum_{i<j} \lk_{Y,\mathbb{Q}}(K_i,K_j)$.  In this case, $\text{sign}(\mathcal{G}_F)+\tfrac{1}{2}e(F,L)$ is the usual signature of $L$ when $Y=\mathbb{S}^3$.

The rational Gordon-Litherland form $\mathcal{G}_F$ yields perhaps the most striking application of the Gordon-Litherland form to date. Answering a question of Fox, Greene gave the following geometric characterization of alternating links. A spanning surface $F$ of a link $L$ is said to be \emph{positive} (\emph{negative}) \emph{definite} if its Gordon-Litherland form is positive (resp. negative) definite.

\begin{theorem}[Greene \cite{greene}] \label{thm_greene} Let $Y$ be a $\mathbb{Z}/2\mathbb{Z}$-homology $3$-sphere. Suppose $L$ is a link in $Y$ having irreducible complement and bounding both a positive definite surface $F_+$ and a negative definite surface $F_-$. Then $Y=\mathbb{S}^3$, $L$ is a non-split alternating link in $\mathbb{S}^3$, and there is a alternating diagram for $L$ whose checkerboard surfaces are isotopic rel boundary to $F_+$, $F_-$. 
\end{theorem}

As a corollary, Greene obtains a proof of one of the Tait conjectures: two connected reduced alternating diagrams of a link have the same number of crossings and the same writhe. This result was originally proved in 1987 by Kauffman \cite{kauffman_87}, Murasugi \cite{murasugi_87}, and Thistlethwaite \cite{thistlethwaite_87}. While the original proof used the Jones polynomial, Greene exploited the Gordon-Litherland pairing to give a more conceptual and geometric argument.

\subsection{Knots in thickened surfaces} 
\begin{figure}[htb]
\begin{tabular}{|ccc|} \hline & & \\
\begin{tabular}{c}
\def\svgwidth{1.5in}
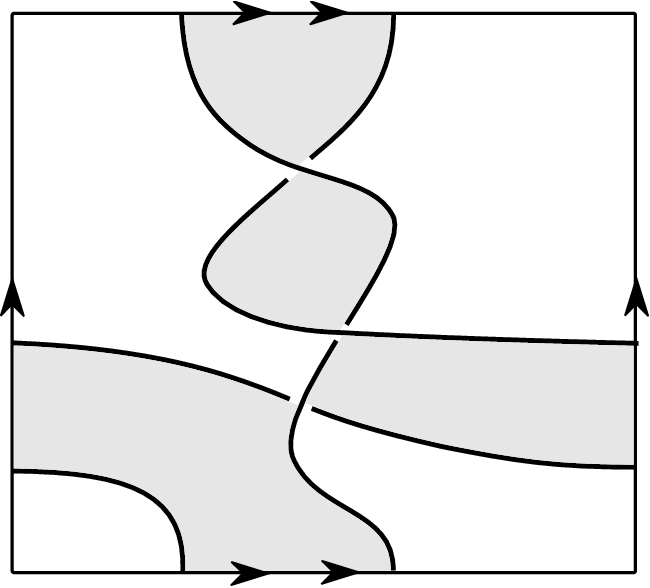
\end{tabular} & &
\begin{tabular}{c}
\def\svgwidth{1.5in}
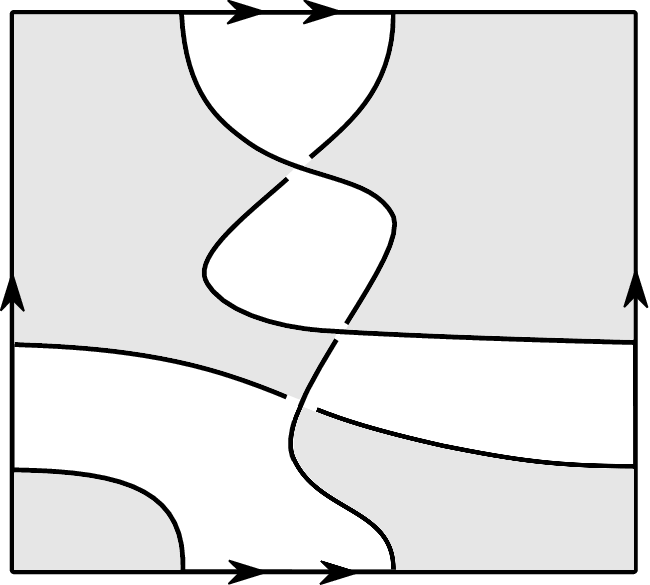  
\end{tabular}
\\ & & \\ \hline
\end{tabular}
\caption{Two checkerboard surfaces of a knot in the thickened torus.} \label{fig_3_7}
\end{figure}

Let $\Sigma$ be a closed oriented surface of genus $g$. A knot $K \subset \Sigma \times [0,1]$ is said to be \emph{checkerboard colorable} if $[K]=0 \in H_1(\Sigma \times [0,1]; \mathbb{Z}/2\mathbb{Z})$. In this case, there is a diagram $D$ of $K$ on $\Sigma$ that is checkerboard colorable in the usual sense. A checkerboard coloring $\xi$ yields two checkerboard surfaces of $K$, $F_{\xi}$ and $F_{\xi}'$ (see Figure \ref{fig_3_7}). The surfaces are said to be \emph{chomatic duals}. For classical knots, chromatic duals are necessarily $S^*$-equivalent. However, this is no longer the case when $g>0$. If $g>0$, every checkerboard colorable knot has exactly two $S^*$-equivalence classes and $F_{\xi},F_{\xi'}$ are never $S^*$-equivalent (see \cite{bck}, Section 1.4). To get from one $S^*$-equivalence class to the other, connect $F_{\xi}$ to a parallel copy of $\Sigma \times \{1\}$ with a compressible tube. The new surface is $S^*$-equivalent to the chromatic dual $F_{\xi'}$ of $F_{\xi}$. 

Since $H_1(\Sigma;\mathbb{Q}) \ne 0$ when $g \ge 1$, the same trick as in Section \ref{sec_z_2z} cannot be used to define the Gordon-Litherland form for knots in thickened surfaces. However, there is an asymmetric linking pairing for such knots due to Cimasoni and Turaev \cite{cimasoni_turaev}. If $K_1,K_2$ are two knots in $\Sigma \times [0,1]$, the linking number $\lk_{\Sigma}(K_1,K_2)$ is defined as the sum of all crossings signs where $J$ overcrosses $K$. If $p:\Sigma \times [0,1] \to \Sigma$ is projection onto the first factor, we have the following identity:
\[
\lk_{\Sigma}(J,K)-\lk_{\Sigma}(K,J)=p_*([J]) \bullet p_*([K]),
\] 
Since the linking number is not symmetric, an additional correction term must be added to the Gordon-Litherland form. If $F$ is any spanning surface for a knot $K \subset \Sigma \times [0,1]$,  $\mathcal{G}_F:H_1(F;\mathbb{Z}) \times H_1(F;\mathbb{Z}) \to \mathbb{Z}$ is defined by (see \cite{bck}, Section 1.2):
\[
\mathcal{G}_F(\alpha_1,\alpha_2)=\lk_{\Sigma}(\tau(\alpha_1),\alpha_2)-p_*(\alpha_1)\bullet p_*(\alpha_2)
\]
With the correction term added, $\mathcal{G}_F$ is now symmetric bilinear form. Define $\sigma(K,F)=\text{sign}(\mathcal{G}_F)+\tfrac{1}{2}e(F)$, where the Euler number is again defined as in Section \ref{sec_gl_pairing}. The signature $\sigma(K,F)$ depends only on the $S^*$-equivalence class of $F$.  Furthermore, in analogy with the classical case, $\mathcal{G}_F$ can be described as a certain relative intersection from of a twofold branched cover along $\widehat{F}$ (see \cite{bck}, Section 6).

A version of Theorem \ref{thm_gl} also holds for knots in thickened surfaces. If $F$ is a Seifert surface of a $\mathbb{Z}$-homologically trivial knot $K$ in $\Sigma \times [0,1]$, $\mathcal{G}_F$ coincides with the symmetrized Seifert form of $F$ (see \cite{bcg}). If $F_{\xi}$ is a checkerboard surface for some coloring $\xi$, one can define the Goeritz matrix $G_{F_{\xi}}$ and the error correction term $\mu(F_{\xi})$ as in Section \ref{sec_what}. Im-Lee-Lee \cite{ILL} defined a signature by the formula $\text{sign}(G_{F_{\xi}})-\mu(F_{\xi})$ and proved that the set of integers: 
\[
\{\text{sign}(G_{F_{\xi}})-\mu(F_{\xi}),\text{sign}(G_{F_{\xi}'})-\mu(F_{\xi'})\},
\]
is an invariant of $K$. In \cite{bck}, it was proved that $\sigma(K,F_{\xi})=G_{F_{\xi'}}-\mu(F_{\xi}')$ and $\sigma(K,F_{\xi'})=G_{F_{\xi}}-\mu(F_{\xi})$. In other words, the Gordon-Litherland signature $\sigma(K,F_{\xi})$ of $F_{\xi}$ equals the Im-Lee-Lee signature of its chromatic dual $F_{\xi'}$ and vice versa.

\bibliographystyle{amsplain}
\bibliography{gl_pairing}

\end{document}

%% file: seven_six_disc_band.eps_tex
\begingroup%
  \makeatletter%
  \providecommand\color[2][]{%
    \errmessage{(Inkscape) Color is used for the text in Inkscape, but the package 'color.sty' is not loaded}%
    \renewcommand\color[2][]{}%
  }%
  \providecommand\transparent[1]{%
    \errmessage{(Inkscape) Transparency is used (non-zero) for the text in Inkscape, but the package 'transparent.sty' is not loaded}%
    \renewcommand\transparent[1]{}%
  }%
  \providecommand\rotatebox[2]{#2}%
  \newcommand*\fsize{\dimexpr\f@size pt\relax}%
  \newcommand*\lineheight[1]{\fontsize{\fsize}{#1\fsize}\selectfont}%
  \ifx\svgwidth\undefined%
    \setlength{\unitlength}{364.40056629bp}%
    \ifx\svgscale\undefined%
      \relax%
    \else%
      \setlength{\unitlength}{\unitlength * \real{\svgscale}}%
    \fi%
  \else%
    \setlength{\unitlength}{\svgwidth}%
  \fi%
  \global\let\svgwidth\undefined%
  \global\let\svgscale\undefined%
  \makeatother%
  \begin{picture}(1,0.69548924)%
    \lineheight{1}%
    \setlength\tabcolsep{0pt}%
    \put(0,0){\includegraphics[width=\unitlength]{seven_six_disc_band.eps}}%
    \put(0.16884816,0.12342346){\color[rgb]{0,0,0}\makebox(0,0)[lt]{\lineheight{40.54999924}\smash{\begin{tabular}[t]{l}$\alpha_1$\end{tabular}}}}%
    \put(0.38052696,0.1218061){\color[rgb]{0,0,0}\makebox(0,0)[lt]{\lineheight{40.54999924}\smash{\begin{tabular}[t]{l}$\alpha_2$\end{tabular}}}}%
    \put(0.60037568,0.1218061){\color[rgb]{0,0,0}\makebox(0,0)[lt]{\lineheight{40.54999924}\smash{\begin{tabular}[t]{l}$\alpha_3$\end{tabular}}}}%
    \put(0.8121416,0.12503923){\color[rgb]{0,0,0}\makebox(0,0)[lt]{\lineheight{40.54999924}\smash{\begin{tabular}[t]{l}$\alpha_4$\end{tabular}}}}%
  \end{picture}%
\endgroup%

%% file: seven_six_check.eps_tex
\begingroup%
  \makeatletter%
  \providecommand\color[2][]{%
    \errmessage{(Inkscape) Color is used for the text in Inkscape, but the package 'color.sty' is not loaded}%
    \renewcommand\color[2][]{}%
  }%
  \providecommand\transparent[1]{%
    \errmessage{(Inkscape) Transparency is used (non-zero) for the text in Inkscape, but the package 'transparent.sty' is not loaded}%
    \renewcommand\transparent[1]{}%
  }%
  \providecommand\rotatebox[2]{#2}%
  \newcommand*\fsize{\dimexpr\f@size pt\relax}%
  \newcommand*\lineheight[1]{\fontsize{\fsize}{#1\fsize}\selectfont}%
  \ifx\svgwidth\undefined%
    \setlength{\unitlength}{224.77925977bp}%
    \ifx\svgscale\undefined%
      \relax%
    \else%
      \setlength{\unitlength}{\unitlength * \real{\svgscale}}%
    \fi%
  \else%
    \setlength{\unitlength}{\svgwidth}%
  \fi%
  \global\let\svgwidth\undefined%
  \global\let\svgscale\undefined%
  \makeatother%
  \begin{picture}(1,1.2214557)%
    \lineheight{1}%
    \setlength\tabcolsep{0pt}%
    \put(0,0){\includegraphics[width=\unitlength]{seven_six_check.eps}}%
    \put(0.84014839,1.09961701){\color[rgb]{0,0,0}\makebox(0,0)[lt]{\lineheight{40.54999924}\smash{\begin{tabular}[t]{l}$X_0$\end{tabular}}}}%
    \put(0.50001549,0.8753517){\color[rgb]{0,0,0}\makebox(0,0)[lt]{\lineheight{40.54999924}\smash{\begin{tabular}[t]{l}$X_1$\end{tabular}}}}%
    \put(0.54518637,0.37359613){\color[rgb]{0,0,0}\makebox(0,0)[lt]{\lineheight{40.54999924}\smash{\begin{tabular}[t]{l}$X_2$\end{tabular}}}}%
    \put(0.12691773,0.38619163){\color[rgb]{0,0,0}\makebox(0,0)[lt]{\lineheight{40.54999924}\smash{\begin{tabular}[t]{l}$X_3$\end{tabular}}}}%
  \end{picture}%
\endgroup%

%% file: seven_six_spanning.eps_tex
\begingroup%
  \makeatletter%
  \providecommand\color[2][]{%
    \errmessage{(Inkscape) Color is used for the text in Inkscape, but the package 'color.sty' is not loaded}%
    \renewcommand\color[2][]{}%
  }%
  \providecommand\transparent[1]{%
    \errmessage{(Inkscape) Transparency is used (non-zero) for the text in Inkscape, but the package 'transparent.sty' is not loaded}%
    \renewcommand\transparent[1]{}%
  }%
  \providecommand\rotatebox[2]{#2}%
  \newcommand*\fsize{\dimexpr\f@size pt\relax}%
  \newcommand*\lineheight[1]{\fontsize{\fsize}{#1\fsize}\selectfont}%
  \ifx\svgwidth\undefined%
    \setlength{\unitlength}{364.40056629bp}%
    \ifx\svgscale\undefined%
      \relax%
    \else%
      \setlength{\unitlength}{\unitlength * \real{\svgscale}}%
    \fi%
  \else%
    \setlength{\unitlength}{\svgwidth}%
  \fi%
  \global\let\svgwidth\undefined%
  \global\let\svgscale\undefined%
  \makeatother%
  \begin{picture}(1,0.69548924)%
    \lineheight{1}%
    \setlength\tabcolsep{0pt}%
    \put(0,0){\includegraphics[width=\unitlength]{seven_six_spanning.eps}}%
    \put(0.17087615,0.1267782){\color[rgb]{0,0,0}\makebox(0,0)[lt]{\lineheight{40.54999924}\smash{\begin{tabular}[t]{l}$\alpha_1$\end{tabular}}}}%
    \put(0.81145359,0.12849045){\color[rgb]{0,0,0}\makebox(0,0)[lt]{\lineheight{40.54999924}\smash{\begin{tabular}[t]{l}$\alpha_3$\end{tabular}}}}%
    \put(0.4962156,0.10379638){\color[rgb]{0,0,0}\makebox(0,0)[lt]{\lineheight{40.54999924}\smash{\begin{tabular}[t]{l}$\alpha_2$\end{tabular}}}}%
  \end{picture}%
\endgroup%

%% file: eta_eq_plus_1.eps_tex
\begingroup%
  \makeatletter%
  \providecommand\color[2][]{%
    \errmessage{(Inkscape) Color is used for the text in Inkscape, but the package 'color.sty' is not loaded}%
    \renewcommand\color[2][]{}%
  }%
  \providecommand\transparent[1]{%
    \errmessage{(Inkscape) Transparency is used (non-zero) for the text in Inkscape, but the package 'transparent.sty' is not loaded}%
    \renewcommand\transparent[1]{}%
  }%
  \providecommand\rotatebox[2]{#2}%
  \newcommand*\fsize{\dimexpr\f@size pt\relax}%
  \newcommand*\lineheight[1]{\fontsize{\fsize}{#1\fsize}\selectfont}%
  \ifx\svgwidth\undefined%
    \setlength{\unitlength}{82.37669307bp}%
    \ifx\svgscale\undefined%
      \relax%
    \else%
      \setlength{\unitlength}{\unitlength * \real{\svgscale}}%
    \fi%
  \else%
    \setlength{\unitlength}{\svgwidth}%
  \fi%
  \global\let\svgwidth\undefined%
  \global\let\svgscale\undefined%
  \makeatother%
  \begin{picture}(1,1)%
    \lineheight{1}%
    \setlength\tabcolsep{0pt}%
    \put(0,0){\includegraphics[width=\unitlength]{eta_eq_plus_1.eps}}%
  \end{picture}%
\endgroup%

%% file: eta_eq_minus_1.eps_tex
\begingroup%
  \makeatletter%
  \providecommand\color[2][]{%
    \errmessage{(Inkscape) Color is used for the text in Inkscape, but the package 'color.sty' is not loaded}%
    \renewcommand\color[2][]{}%
  }%
  \providecommand\transparent[1]{%
    \errmessage{(Inkscape) Transparency is used (non-zero) for the text in Inkscape, but the package 'transparent.sty' is not loaded}%
    \renewcommand\transparent[1]{}%
  }%
  \providecommand\rotatebox[2]{#2}%
  \newcommand*\fsize{\dimexpr\f@size pt\relax}%
  \newcommand*\lineheight[1]{\fontsize{\fsize}{#1\fsize}\selectfont}%
  \ifx\svgwidth\undefined%
    \setlength{\unitlength}{82.37669307bp}%
    \ifx\svgscale\undefined%
      \relax%
    \else%
      \setlength{\unitlength}{\unitlength * \real{\svgscale}}%
    \fi%
  \else%
    \setlength{\unitlength}{\svgwidth}%
  \fi%
  \global\let\svgwidth\undefined%
  \global\let\svgscale\undefined%
  \makeatother%
  \begin{picture}(1,1)%
    \lineheight{1}%
    \setlength\tabcolsep{0pt}%
    \put(0,0){\includegraphics[width=\unitlength]{eta_eq_minus_1.eps}}%
  \end{picture}%
\endgroup%

%% file: type_I.eps_tex
\begingroup%
  \makeatletter%
  \providecommand\color[2][]{%
    \errmessage{(Inkscape) Color is used for the text in Inkscape, but the package 'color.sty' is not loaded}%
    \renewcommand\color[2][]{}%
  }%
  \providecommand\transparent[1]{%
    \errmessage{(Inkscape) Transparency is used (non-zero) for the text in Inkscape, but the package 'transparent.sty' is not loaded}%
    \renewcommand\transparent[1]{}%
  }%
  \providecommand\rotatebox[2]{#2}%
  \newcommand*\fsize{\dimexpr\f@size pt\relax}%
  \newcommand*\lineheight[1]{\fontsize{\fsize}{#1\fsize}\selectfont}%
  \ifx\svgwidth\undefined%
    \setlength{\unitlength}{84.04264793bp}%
    \ifx\svgscale\undefined%
      \relax%
    \else%
      \setlength{\unitlength}{\unitlength * \real{\svgscale}}%
    \fi%
  \else%
    \setlength{\unitlength}{\svgwidth}%
  \fi%
  \global\let\svgwidth\undefined%
  \global\let\svgscale\undefined%
  \makeatother%
  \begin{picture}(1,1.0198227)%
    \lineheight{1}%
    \setlength\tabcolsep{0pt}%
    \put(0,0){\includegraphics[width=\unitlength]{type_I.eps}}%
  \end{picture}%
\endgroup%

%% file: type_II.eps_tex
\begingroup%
  \makeatletter%
  \providecommand\color[2][]{%
    \errmessage{(Inkscape) Color is used for the text in Inkscape, but the package 'color.sty' is not loaded}%
    \renewcommand\color[2][]{}%
  }%
  \providecommand\transparent[1]{%
    \errmessage{(Inkscape) Transparency is used (non-zero) for the text in Inkscape, but the package 'transparent.sty' is not loaded}%
    \renewcommand\transparent[1]{}%
  }%
  \providecommand\rotatebox[2]{#2}%
  \newcommand*\fsize{\dimexpr\f@size pt\relax}%
  \newcommand*\lineheight[1]{\fontsize{\fsize}{#1\fsize}\selectfont}%
  \ifx\svgwidth\undefined%
    \setlength{\unitlength}{85.70859969bp}%
    \ifx\svgscale\undefined%
      \relax%
    \else%
      \setlength{\unitlength}{\unitlength * \real{\svgscale}}%
    \fi%
  \else%
    \setlength{\unitlength}{\svgwidth}%
  \fi%
  \global\let\svgwidth\undefined%
  \global\let\svgscale\undefined%
  \makeatother%
  \begin{picture}(1,0.98056261)%
    \lineheight{1}%
    \setlength\tabcolsep{0pt}%
    \put(0,0){\includegraphics[width=\unitlength]{type_II.eps}}%
  \end{picture}%
\endgroup%

%% file: goeritz_basis.eps_tex
\begingroup%
  \makeatletter%
  \providecommand\color[2][]{%
    \errmessage{(Inkscape) Color is used for the text in Inkscape, but the package 'color.sty' is not loaded}%
    \renewcommand\color[2][]{}%
  }%
  \providecommand\transparent[1]{%
    \errmessage{(Inkscape) Transparency is used (non-zero) for the text in Inkscape, but the package 'transparent.sty' is not loaded}%
    \renewcommand\transparent[1]{}%
  }%
  \providecommand\rotatebox[2]{#2}%
  \newcommand*\fsize{\dimexpr\f@size pt\relax}%
  \newcommand*\lineheight[1]{\fontsize{\fsize}{#1\fsize}\selectfont}%
  \ifx\svgwidth\undefined%
    \setlength{\unitlength}{239.0779125bp}%
    \ifx\svgscale\undefined%
      \relax%
    \else%
      \setlength{\unitlength}{\unitlength * \real{\svgscale}}%
    \fi%
  \else%
    \setlength{\unitlength}{\svgwidth}%
  \fi%
  \global\let\svgwidth\undefined%
  \global\let\svgscale\undefined%
  \makeatother%
  \begin{picture}(1,0.81517004)%
    \lineheight{1}%
    \setlength\tabcolsep{0pt}%
    \put(0,0){\includegraphics[width=\unitlength]{goeritz_basis.eps}}%
    \put(0.35844155,0.39341271){\color[rgb]{0,0,0}\makebox(0,0)[lt]{\lineheight{40.54999924}\smash{\begin{tabular}[t]{l}$X_i$\end{tabular}}}}%
    \put(0.77099379,0.39523007){\color[rgb]{0,0,0}\makebox(0,0)[lt]{\lineheight{40.54999924}\smash{\begin{tabular}[t]{l}$\alpha_i$\end{tabular}}}}%
  \end{picture}%
\endgroup%

%% file: goeritz_entries.eps_tex
\begingroup%
  \makeatletter%
  \providecommand\color[2][]{%
    \errmessage{(Inkscape) Color is used for the text in Inkscape, but the package 'color.sty' is not loaded}%
    \renewcommand\color[2][]{}%
  }%
  \providecommand\transparent[1]{%
    \errmessage{(Inkscape) Transparency is used (non-zero) for the text in Inkscape, but the package 'transparent.sty' is not loaded}%
    \renewcommand\transparent[1]{}%
  }%
  \providecommand\rotatebox[2]{#2}%
  \newcommand*\fsize{\dimexpr\f@size pt\relax}%
  \newcommand*\lineheight[1]{\fontsize{\fsize}{#1\fsize}\selectfont}%
  \ifx\svgwidth\undefined%
    \setlength{\unitlength}{180.96508953bp}%
    \ifx\svgscale\undefined%
      \relax%
    \else%
      \setlength{\unitlength}{\unitlength * \real{\svgscale}}%
    \fi%
  \else%
    \setlength{\unitlength}{\svgwidth}%
  \fi%
  \global\let\svgwidth\undefined%
  \global\let\svgscale\undefined%
  \makeatother%
  \begin{picture}(1,0.94089723)%
    \lineheight{1}%
    \setlength\tabcolsep{0pt}%
    \put(0,0){\includegraphics[width=\unitlength]{goeritz_entries.eps}}%
    \put(0.59285273,0.45597282){\color[rgb]{0,0,0}\makebox(0,0)[lt]{\lineheight{40.54999924}\smash{\begin{tabular}[t]{l}$X_j$\end{tabular}}}}%
    \put(0.10069667,0.45373571){\color[rgb]{0,0,0}\makebox(0,0)[lt]{\lineheight{40.54999924}\smash{\begin{tabular}[t]{l}$X_i$\end{tabular}}}}%
    \put(0.00002838,0.03316599){\color[rgb]{0,0,0}\makebox(0,0)[lt]{\lineheight{40.54999924}\smash{\begin{tabular}[t]{l}$\alpha_i$\end{tabular}}}}%
    \put(0.63088296,0.01303231){\color[rgb]{0,0,0}\makebox(0,0)[lt]{\lineheight{40.54999924}\smash{\begin{tabular}[t]{l} \textcolor{red}{$\alpha_j$} \end{tabular}}}}%
    \put(0.54811126,0.89891322){\color[rgb]{0,0,0}\makebox(0,0)[lt]{\lineheight{40.54999924}\smash{\begin{tabular}[t]{l}\textcolor[RGB]{13,170,13}{$\tau(\alpha_j)$}\end{tabular}}}}%
  \end{picture}%
\endgroup%

%% file: goeritz_entries_II.eps_tex
\begingroup%
  \makeatletter%
  \providecommand\color[2][]{%
    \errmessage{(Inkscape) Color is used for the text in Inkscape, but the package 'color.sty' is not loaded}%
    \renewcommand\color[2][]{}%
  }%
  \providecommand\transparent[1]{%
    \errmessage{(Inkscape) Transparency is used (non-zero) for the text in Inkscape, but the package 'transparent.sty' is not loaded}%
    \renewcommand\transparent[1]{}%
  }%
  \providecommand\rotatebox[2]{#2}%
  \newcommand*\fsize{\dimexpr\f@size pt\relax}%
  \newcommand*\lineheight[1]{\fontsize{\fsize}{#1\fsize}\selectfont}%
  \ifx\svgwidth\undefined%
    \setlength{\unitlength}{173.47246812bp}%
    \ifx\svgscale\undefined%
      \relax%
    \else%
      \setlength{\unitlength}{\unitlength * \real{\svgscale}}%
    \fi%
  \else%
    \setlength{\unitlength}{\svgwidth}%
  \fi%
  \global\let\svgwidth\undefined%
  \global\let\svgscale\undefined%
  \makeatother%
  \begin{picture}(1,0.97644224)%
    \lineheight{1}%
    \setlength\tabcolsep{0pt}%
    \put(0,0){\includegraphics[width=\unitlength]{goeritz_entries_II.eps}}%
    \put(0.54902927,0.487391){\color[rgb]{0,0,0}\makebox(0,0)[lt]{\lineheight{40.54999924}\smash{\begin{tabular}[t]{l}$X_j$\end{tabular}}}}%
    \put(0.08855303,0.47026585){\color[rgb]{0,0,0}\makebox(0,0)[lt]{\lineheight{40.54999924}\smash{\begin{tabular}[t]{l}$X_i$\end{tabular}}}}%
    \put(0.02766361,0.01549802){\color[rgb]{0,0,0}\makebox(0,0)[lt]{\lineheight{40.54999924}\smash{\begin{tabular}[t]{l}$\alpha_i$\end{tabular}}}}%
    \put(0.62704385,0.0135952){\color[rgb]{0,0,0}\makebox(0,0)[lt]{\lineheight{40.54999924}\smash{\begin{tabular}[t]{l} \textcolor{red}{$\alpha_j$}\end{tabular}}}}%
    \put(0.51477897,0.93264486){\color[rgb]{0,0,0}\makebox(0,0)[lt]{\lineheight{40.54999924}\smash{\begin{tabular}[t]{l} \textcolor[RGB]{13,170,13}{$\tau(\alpha_j)$}\end{tabular}}}}%
  \end{picture}%
\endgroup%

%% file: eta_eq_1_type_eq_1.eps_tex
\begingroup%
  \makeatletter%
  \providecommand\color[2][]{%
    \errmessage{(Inkscape) Color is used for the text in Inkscape, but the package 'color.sty' is not loaded}%
    \renewcommand\color[2][]{}%
  }%
  \providecommand\transparent[1]{%
    \errmessage{(Inkscape) Transparency is used (non-zero) for the text in Inkscape, but the package 'transparent.sty' is not loaded}%
    \renewcommand\transparent[1]{}%
  }%
  \providecommand\rotatebox[2]{#2}%
  \newcommand*\fsize{\dimexpr\f@size pt\relax}%
  \newcommand*\lineheight[1]{\fontsize{\fsize}{#1\fsize}\selectfont}%
  \ifx\svgwidth\undefined%
    \setlength{\unitlength}{100.89689673bp}%
    \ifx\svgscale\undefined%
      \relax%
    \else%
      \setlength{\unitlength}{\unitlength * \real{\svgscale}}%
    \fi%
  \else%
    \setlength{\unitlength}{\svgwidth}%
  \fi%
  \global\let\svgwidth\undefined%
  \global\let\svgscale\undefined%
  \makeatother%
  \begin{picture}(1,0.9127645)%
    \lineheight{1}%
    \setlength\tabcolsep{0pt}%
    \put(0,0){\includegraphics[width=\unitlength]{eta_eq_1_type_eq_1.eps}}%
    \put(0.72439417,0.41239662){\color[rgb]{0,0,0}\makebox(0,0)[lt]{\lineheight{40.54999924}\smash{\begin{tabular}[t]{l}$K'$\end{tabular}}}}%
    \put(0.5983693,0.83746348){\color[rgb]{0,0,0}\makebox(0,0)[lt]{\lineheight{40.54999924}\smash{\begin{tabular}[t]{l}$K$\end{tabular}}}}%
  \end{picture}%
\endgroup%

%% file: eta_eq_neg_1_type_eq_1.eps_tex
\begingroup%
  \makeatletter%
  \providecommand\color[2][]{%
    \errmessage{(Inkscape) Color is used for the text in Inkscape, but the package 'color.sty' is not loaded}%
    \renewcommand\color[2][]{}%
  }%
  \providecommand\transparent[1]{%
    \errmessage{(Inkscape) Transparency is used (non-zero) for the text in Inkscape, but the package 'transparent.sty' is not loaded}%
    \renewcommand\transparent[1]{}%
  }%
  \providecommand\rotatebox[2]{#2}%
  \newcommand*\fsize{\dimexpr\f@size pt\relax}%
  \newcommand*\lineheight[1]{\fontsize{\fsize}{#1\fsize}\selectfont}%
  \ifx\svgwidth\undefined%
    \setlength{\unitlength}{100.89689778bp}%
    \ifx\svgscale\undefined%
      \relax%
    \else%
      \setlength{\unitlength}{\unitlength * \real{\svgscale}}%
    \fi%
  \else%
    \setlength{\unitlength}{\svgwidth}%
  \fi%
  \global\let\svgwidth\undefined%
  \global\let\svgscale\undefined%
  \makeatother%
  \begin{picture}(1,0.9170365)%
    \lineheight{1}%
    \setlength\tabcolsep{0pt}%
    \put(0,0){\includegraphics[width=\unitlength]{eta_eq_neg_1_type_eq_1.eps}}%
    \put(0.73080218,0.41666858){\color[rgb]{0,0,0}\makebox(0,0)[lt]{\lineheight{40.54999924}\smash{\begin{tabular}[t]{l}$K'$\end{tabular}}}}%
    \put(0.60477743,0.84173548){\color[rgb]{0,0,0}\makebox(0,0)[lt]{\lineheight{40.54999924}\smash{\begin{tabular}[t]{l}$K$\end{tabular}}}}%
  \end{picture}%
\endgroup%

%% file: eta_eq_1_type_eq_2.eps_tex
\begingroup%
  \makeatletter%
  \providecommand\color[2][]{%
    \errmessage{(Inkscape) Color is used for the text in Inkscape, but the package 'color.sty' is not loaded}%
    \renewcommand\color[2][]{}%
  }%
  \providecommand\transparent[1]{%
    \errmessage{(Inkscape) Transparency is used (non-zero) for the text in Inkscape, but the package 'transparent.sty' is not loaded}%
    \renewcommand\transparent[1]{}%
  }%
  \providecommand\rotatebox[2]{#2}%
  \newcommand*\fsize{\dimexpr\f@size pt\relax}%
  \newcommand*\lineheight[1]{\fontsize{\fsize}{#1\fsize}\selectfont}%
  \ifx\svgwidth\undefined%
    \setlength{\unitlength}{102.2602085bp}%
    \ifx\svgscale\undefined%
      \relax%
    \else%
      \setlength{\unitlength}{\unitlength * \real{\svgscale}}%
    \fi%
  \else%
    \setlength{\unitlength}{\svgwidth}%
  \fi%
  \global\let\svgwidth\undefined%
  \global\let\svgscale\undefined%
  \makeatother%
  \begin{picture}(1,0.88642349)%
    \lineheight{1}%
    \setlength\tabcolsep{0pt}%
    \put(0,0){\includegraphics[width=\unitlength]{eta_eq_1_type_eq_2.eps}}%
    \put(0.77654177,0.39272636){\color[rgb]{0,0,0}\makebox(0,0)[lt]{\lineheight{40.54999924}\smash{\begin{tabular}[t]{l}$K'$\end{tabular}}}}%
    \put(0.6521972,0.81212636){\color[rgb]{0,0,0}\makebox(0,0)[lt]{\lineheight{40.54999924}\smash{\begin{tabular}[t]{l}$K$\end{tabular}}}}%
  \end{picture}%
\endgroup%

%% file: eta_eq_neg_1_type_eq_2.eps_tex
\begingroup%
  \makeatletter%
  \providecommand\color[2][]{%
    \errmessage{(Inkscape) Color is used for the text in Inkscape, but the package 'color.sty' is not loaded}%
    \renewcommand\color[2][]{}%
  }%
  \providecommand\transparent[1]{%
    \errmessage{(Inkscape) Transparency is used (non-zero) for the text in Inkscape, but the package 'transparent.sty' is not loaded}%
    \renewcommand\transparent[1]{}%
  }%
  \providecommand\rotatebox[2]{#2}%
  \newcommand*\fsize{\dimexpr\f@size pt\relax}%
  \newcommand*\lineheight[1]{\fontsize{\fsize}{#1\fsize}\selectfont}%
  \ifx\svgwidth\undefined%
    \setlength{\unitlength}{102.26020901bp}%
    \ifx\svgscale\undefined%
      \relax%
    \else%
      \setlength{\unitlength}{\unitlength * \real{\svgscale}}%
    \fi%
  \else%
    \setlength{\unitlength}{\svgwidth}%
  \fi%
  \global\let\svgwidth\undefined%
  \global\let\svgscale\undefined%
  \makeatother%
  \begin{picture}(1,0.89425669)%
    \lineheight{1}%
    \setlength\tabcolsep{0pt}%
    \put(0,0){\includegraphics[width=\unitlength]{eta_eq_neg_1_type_eq_2.eps}}%
    \put(0.77361733,0.4005596){\color[rgb]{0,0,0}\makebox(0,0)[lt]{\lineheight{40.54999924}\smash{\begin{tabular}[t]{l}$K'$\end{tabular}}}}%
    \put(0.64927266,0.81995957){\color[rgb]{0,0,0}\makebox(0,0)[lt]{\lineheight{40.54999924}\smash{\begin{tabular}[t]{l}$K$\end{tabular}}}}%
  \end{picture}%
\endgroup%

%% file: pos_curl.eps_tex
\begingroup%
  \makeatletter%
  \providecommand\color[2][]{%
    \errmessage{(Inkscape) Color is used for the text in Inkscape, but the package 'color.sty' is not loaded}%
    \renewcommand\color[2][]{}%
  }%
  \providecommand\transparent[1]{%
    \errmessage{(Inkscape) Transparency is used (non-zero) for the text in Inkscape, but the package 'transparent.sty' is not loaded}%
    \renewcommand\transparent[1]{}%
  }%
  \providecommand\rotatebox[2]{#2}%
  \newcommand*\fsize{\dimexpr\f@size pt\relax}%
  \newcommand*\lineheight[1]{\fontsize{\fsize}{#1\fsize}\selectfont}%
  \ifx\svgwidth\undefined%
    \setlength{\unitlength}{68.26445788bp}%
    \ifx\svgscale\undefined%
      \relax%
    \else%
      \setlength{\unitlength}{\unitlength * \real{\svgscale}}%
    \fi%
  \else%
    \setlength{\unitlength}{\svgwidth}%
  \fi%
  \global\let\svgwidth\undefined%
  \global\let\svgscale\undefined%
  \makeatother%
  \begin{picture}(1,1.43869103)%
    \lineheight{1}%
    \setlength\tabcolsep{0pt}%
    \put(0,0){\includegraphics[width=\unitlength]{pos_curl.eps}}%
  \end{picture}%
\endgroup%

%% file: pos_full_twist.eps_tex
\begingroup%
  \makeatletter%
  \providecommand\color[2][]{%
    \errmessage{(Inkscape) Color is used for the text in Inkscape, but the package 'color.sty' is not loaded}%
    \renewcommand\color[2][]{}%
  }%
  \providecommand\transparent[1]{%
    \errmessage{(Inkscape) Transparency is used (non-zero) for the text in Inkscape, but the package 'transparent.sty' is not loaded}%
    \renewcommand\transparent[1]{}%
  }%
  \providecommand\rotatebox[2]{#2}%
  \newcommand*\fsize{\dimexpr\f@size pt\relax}%
  \newcommand*\lineheight[1]{\fontsize{\fsize}{#1\fsize}\selectfont}%
  \ifx\svgwidth\undefined%
    \setlength{\unitlength}{104.25902893bp}%
    \ifx\svgscale\undefined%
      \relax%
    \else%
      \setlength{\unitlength}{\unitlength * \real{\svgscale}}%
    \fi%
  \else%
    \setlength{\unitlength}{\svgwidth}%
  \fi%
  \global\let\svgwidth\undefined%
  \global\let\svgscale\undefined%
  \makeatother%
  \begin{picture}(1,0.2257413)%
    \lineheight{1}%
    \setlength\tabcolsep{0pt}%
    \put(0,0){\includegraphics[width=\unitlength]{pos_full_twist.eps}}%
  \end{picture}%
\endgroup%

%% file: neg_curl.eps_tex
\begingroup%
  \makeatletter%
  \providecommand\color[2][]{%
    \errmessage{(Inkscape) Color is used for the text in Inkscape, but the package 'color.sty' is not loaded}%
    \renewcommand\color[2][]{}%
  }%
  \providecommand\transparent[1]{%
    \errmessage{(Inkscape) Transparency is used (non-zero) for the text in Inkscape, but the package 'transparent.sty' is not loaded}%
    \renewcommand\transparent[1]{}%
  }%
  \providecommand\rotatebox[2]{#2}%
  \newcommand*\fsize{\dimexpr\f@size pt\relax}%
  \newcommand*\lineheight[1]{\fontsize{\fsize}{#1\fsize}\selectfont}%
  \ifx\svgwidth\undefined%
    \setlength{\unitlength}{68.2644593bp}%
    \ifx\svgscale\undefined%
      \relax%
    \else%
      \setlength{\unitlength}{\unitlength * \real{\svgscale}}%
    \fi%
  \else%
    \setlength{\unitlength}{\svgwidth}%
  \fi%
  \global\let\svgwidth\undefined%
  \global\let\svgscale\undefined%
  \makeatother%
  \begin{picture}(1,1.438691)%
    \lineheight{1}%
    \setlength\tabcolsep{0pt}%
    \put(0,0){\includegraphics[width=\unitlength]{neg_curl.eps}}%
  \end{picture}%
\endgroup%

%% file: neg_full_twist.eps_tex
\begingroup%
  \makeatletter%
  \providecommand\color[2][]{%
    \errmessage{(Inkscape) Color is used for the text in Inkscape, but the package 'color.sty' is not loaded}%
    \renewcommand\color[2][]{}%
  }%
  \providecommand\transparent[1]{%
    \errmessage{(Inkscape) Transparency is used (non-zero) for the text in Inkscape, but the package 'transparent.sty' is not loaded}%
    \renewcommand\transparent[1]{}%
  }%
  \providecommand\rotatebox[2]{#2}%
  \newcommand*\fsize{\dimexpr\f@size pt\relax}%
  \newcommand*\lineheight[1]{\fontsize{\fsize}{#1\fsize}\selectfont}%
  \ifx\svgwidth\undefined%
    \setlength{\unitlength}{104.25902893bp}%
    \ifx\svgscale\undefined%
      \relax%
    \else%
      \setlength{\unitlength}{\unitlength * \real{\svgscale}}%
    \fi%
  \else%
    \setlength{\unitlength}{\svgwidth}%
  \fi%
  \global\let\svgwidth\undefined%
  \global\let\svgscale\undefined%
  \makeatother%
  \begin{picture}(1,0.2257413)%
    \lineheight{1}%
    \setlength\tabcolsep{0pt}%
    \put(0,0){\includegraphics[width=\unitlength]{neg_full_twist.eps}}%
  \end{picture}%
\endgroup%

%% file: seven_six_kirby_I.eps_tex
\begingroup%
  \makeatletter%
  \providecommand\color[2][]{%
    \errmessage{(Inkscape) Color is used for the text in Inkscape, but the package 'color.sty' is not loaded}%
    \renewcommand\color[2][]{}%
  }%
  \providecommand\transparent[1]{%
    \errmessage{(Inkscape) Transparency is used (non-zero) for the text in Inkscape, but the package 'transparent.sty' is not loaded}%
    \renewcommand\transparent[1]{}%
  }%
  \providecommand\rotatebox[2]{#2}%
  \newcommand*\fsize{\dimexpr\f@size pt\relax}%
  \newcommand*\lineheight[1]{\fontsize{\fsize}{#1\fsize}\selectfont}%
  \ifx\svgwidth\undefined%
    \setlength{\unitlength}{433.12858256bp}%
    \ifx\svgscale\undefined%
      \relax%
    \else%
      \setlength{\unitlength}{\unitlength * \real{\svgscale}}%
    \fi%
  \else%
    \setlength{\unitlength}{\svgwidth}%
  \fi%
  \global\let\svgwidth\undefined%
  \global\let\svgscale\undefined%
  \makeatother%
  \begin{picture}(1,0.58908689)%
    \lineheight{1}%
    \setlength\tabcolsep{0pt}%
    \put(0,0){\includegraphics[width=\unitlength]{seven_six_kirby_I.eps}}%
    \put(0.22496501,0.11633685){\color[rgb]{0,0,0}\makebox(0,0)[lt]{\lineheight{40.54999924}\smash{\begin{tabular}[t]{l}$\alpha_1$\end{tabular}}}}%
    \put(0.76389684,0.11920719){\color[rgb]{0,0,0}\makebox(0,0)[lt]{\lineheight{40.54999924}\smash{\begin{tabular}[t]{l}$\alpha_3$\end{tabular}}}}%
    \put(0.49868021,0.09914639){\color[rgb]{0,0,0}\makebox(0,0)[lt]{\lineheight{40.54999924}\smash{\begin{tabular}[t]{l}$\alpha_2$\end{tabular}}}}%
    \put(0.1014101,0.36698396){\color[rgb]{0,0,0}\makebox(0,0)[lt]{\lineheight{40.54999924}\smash{\begin{tabular}[t]{l}$x$\end{tabular}}}}%
    \put(0.04805385,0.44237868){\color[rgb]{0,0,0}\makebox(0,0)[lt]{\lineheight{40.54999924}\smash{\begin{tabular}[t]{l}$y$\end{tabular}}}}%
    \put(-0.00182264,0.39876575){\color[rgb]{0,0,0}\makebox(0,0)[lt]{\lineheight{40.54999924}\smash{\begin{tabular}[t]{l}$z$\end{tabular}}}}%
    \put(0.0106108,0.19294531){\color[rgb]{0,0,0}\rotatebox{43.2}{\makebox(0,0)[lt]{\lineheight{40.54999924}\smash{\begin{tabular}[t]{l}$x$-axis\end{tabular}}}}}%
    \put(0.93737944,0.1376013){\color[rgb]{0,0,0}\makebox(0,0)[lt]{\lineheight{40.54999924}\smash{\begin{tabular}[t]{l}$\pi$\end{tabular}}}}%
  \end{picture}%
\endgroup%

%% file: seven_six_framed_link.eps_tex
\begingroup%
  \makeatletter%
  \providecommand\color[2][]{%
    \errmessage{(Inkscape) Color is used for the text in Inkscape, but the package 'color.sty' is not loaded}%
    \renewcommand\color[2][]{}%
  }%
  \providecommand\transparent[1]{%
    \errmessage{(Inkscape) Transparency is used (non-zero) for the text in Inkscape, but the package 'transparent.sty' is not loaded}%
    \renewcommand\transparent[1]{}%
  }%
  \providecommand\rotatebox[2]{#2}%
  \newcommand*\fsize{\dimexpr\f@size pt\relax}%
  \newcommand*\lineheight[1]{\fontsize{\fsize}{#1\fsize}\selectfont}%
  \ifx\svgwidth\undefined%
    \setlength{\unitlength}{253.6652922bp}%
    \ifx\svgscale\undefined%
      \relax%
    \else%
      \setlength{\unitlength}{\unitlength * \real{\svgscale}}%
    \fi%
  \else%
    \setlength{\unitlength}{\svgwidth}%
  \fi%
  \global\let\svgwidth\undefined%
  \global\let\svgscale\undefined%
  \makeatother%
  \begin{picture}(1,0.50282)%
    \lineheight{1}%
    \setlength\tabcolsep{0pt}%
    \put(0,0){\includegraphics[width=\unitlength]{seven_six_framed_link.eps}}%
    \put(0.20982575,0.46736378){\color[rgb]{0,0,0}\makebox(0,0)[lt]{\lineheight{40.54999924}\smash{\begin{tabular}[t]{l}$3$\end{tabular}}}}%
    \put(0.4860832,0.47286856){\color[rgb]{0,0,0}\makebox(0,0)[lt]{\lineheight{40.54999924}\smash{\begin{tabular}[t]{l}$4$\end{tabular}}}}%
    \put(0.75802647,0.47208388){\color[rgb]{0,0,0}\makebox(0,0)[lt]{\lineheight{40.54999924}\smash{\begin{tabular}[t]{l}$2$\end{tabular}}}}%
  \end{picture}%
\endgroup%

%% file: bubble_gum_left.eps_tex
\begingroup%
  \makeatletter%
  \providecommand\color[2][]{%
    \errmessage{(Inkscape) Color is used for the text in Inkscape, but the package 'color.sty' is not loaded}%
    \renewcommand\color[2][]{}%
  }%
  \providecommand\transparent[1]{%
    \errmessage{(Inkscape) Transparency is used (non-zero) for the text in Inkscape, but the package 'transparent.sty' is not loaded}%
    \renewcommand\transparent[1]{}%
  }%
  \providecommand\rotatebox[2]{#2}%
  \newcommand*\fsize{\dimexpr\f@size pt\relax}%
  \newcommand*\lineheight[1]{\fontsize{\fsize}{#1\fsize}\selectfont}%
  \ifx\svgwidth\undefined%
    \setlength{\unitlength}{222.44092434bp}%
    \ifx\svgscale\undefined%
      \relax%
    \else%
      \setlength{\unitlength}{\unitlength * \real{\svgscale}}%
    \fi%
  \else%
    \setlength{\unitlength}{\svgwidth}%
  \fi%
  \global\let\svgwidth\undefined%
  \global\let\svgscale\undefined%
  \makeatother%
  \begin{picture}(1,0.58003634)%
    \lineheight{1}%
    \setlength\tabcolsep{0pt}%
    \put(0,0){\includegraphics[width=\unitlength]{bubble_gum_left.eps}}%
    \put(0.38771885,0.54588057){\color[rgb]{0,0,0}\makebox(0,0)[lt]{\lineheight{40.54999924}\smash{\begin{tabular}[t]{l}$\mathbb{B}^2 \times \{0\}$\end{tabular}}}}%
    \put(0.16186217,0.022){\color[rgb]{0,0,0}\makebox(0,0)[lt]{\lineheight{40.54999924}\smash{\begin{tabular}[t]{l}$F^1_i$\end{tabular}}}}%
    \put(0.62084354,0.022){\color[rgb]{0,0,0}\makebox(0,0)[lt]{\lineheight{40.54999924}\smash{\begin{tabular}[t]{l}$V^1_i$\end{tabular}}}}%
  \end{picture}%
\endgroup%

%% file: tube_left.eps_tex
\begingroup%
  \makeatletter%
  \providecommand\color[2][]{%
    \errmessage{(Inkscape) Color is used for the text in Inkscape, but the package 'color.sty' is not loaded}%
    \renewcommand\color[2][]{}%
  }%
  \providecommand\transparent[1]{%
    \errmessage{(Inkscape) Transparency is used (non-zero) for the text in Inkscape, but the package 'transparent.sty' is not loaded}%
    \renewcommand\transparent[1]{}%
  }%
  \providecommand\rotatebox[2]{#2}%
  \newcommand*\fsize{\dimexpr\f@size pt\relax}%
  \newcommand*\lineheight[1]{\fontsize{\fsize}{#1\fsize}\selectfont}%
  \ifx\svgwidth\undefined%
    \setlength{\unitlength}{225.65257995bp}%
    \ifx\svgscale\undefined%
      \relax%
    \else%
      \setlength{\unitlength}{\unitlength * \real{\svgscale}}%
    \fi%
  \else%
    \setlength{\unitlength}{\svgwidth}%
  \fi%
  \global\let\svgwidth\undefined%
  \global\let\svgscale\undefined%
  \makeatother%
  \begin{picture}(1,0.1713233)%
    \lineheight{1}%
    \setlength\tabcolsep{0pt}%
    \put(0,0){\includegraphics[width=\unitlength]{tube_left.eps}}%
  \end{picture}%
\endgroup%

%% file: tube_right.eps_tex
\begingroup%
  \makeatletter%
  \providecommand\color[2][]{%
    \errmessage{(Inkscape) Color is used for the text in Inkscape, but the package 'color.sty' is not loaded}%
    \renewcommand\color[2][]{}%
  }%
  \providecommand\transparent[1]{%
    \errmessage{(Inkscape) Transparency is used (non-zero) for the text in Inkscape, but the package 'transparent.sty' is not loaded}%
    \renewcommand\transparent[1]{}%
  }%
  \providecommand\rotatebox[2]{#2}%
  \newcommand*\fsize{\dimexpr\f@size pt\relax}%
  \newcommand*\lineheight[1]{\fontsize{\fsize}{#1\fsize}\selectfont}%
  \ifx\svgwidth\undefined%
    \setlength{\unitlength}{292.35903845bp}%
    \ifx\svgscale\undefined%
      \relax%
    \else%
      \setlength{\unitlength}{\unitlength * \real{\svgscale}}%
    \fi%
  \else%
    \setlength{\unitlength}{\svgwidth}%
  \fi%
  \global\let\svgwidth\undefined%
  \global\let\svgscale\undefined%
  \makeatother%
  \begin{picture}(1,0.32577711)%
    \lineheight{1}%
    \setlength\tabcolsep{0pt}%
    \put(0,0){\includegraphics[width=\unitlength]{tube_right.eps}}%
    \put(0.45962565,0.23438657){\color[rgb]{0,0,0}\makebox(0,0)[lt]{\lineheight{40.54999924}\smash{\begin{tabular}[t]{l}$\alpha_{n+2}$\end{tabular}}}}%
    \put(0.72025173,0.02674493){\color[rgb]{0,0,0}\makebox(0,0)[lt]{\lineheight{40.54999924}\smash{\begin{tabular}[t]{l}$\alpha_{n+1}$\end{tabular}}}}%
  \end{picture}%
\endgroup%

%% file: twist_right_I.eps_tex
\begingroup%
  \makeatletter%
  \providecommand\color[2][]{%
    \errmessage{(Inkscape) Color is used for the text in Inkscape, but the package 'color.sty' is not loaded}%
    \renewcommand\color[2][]{}%
  }%
  \providecommand\transparent[1]{%
    \errmessage{(Inkscape) Transparency is used (non-zero) for the text in Inkscape, but the package 'transparent.sty' is not loaded}%
    \renewcommand\transparent[1]{}%
  }%
  \providecommand\rotatebox[2]{#2}%
  \newcommand*\fsize{\dimexpr\f@size pt\relax}%
  \newcommand*\lineheight[1]{\fontsize{\fsize}{#1\fsize}\selectfont}%
  \ifx\svgwidth\undefined%
    \setlength{\unitlength}{162.44276172bp}%
    \ifx\svgscale\undefined%
      \relax%
    \else%
      \setlength{\unitlength}{\unitlength * \real{\svgscale}}%
    \fi%
  \else%
    \setlength{\unitlength}{\svgwidth}%
  \fi%
  \global\let\svgwidth\undefined%
  \global\let\svgscale\undefined%
  \makeatother%
  \begin{picture}(1,0.96339347)%
    \lineheight{1}%
    \setlength\tabcolsep{0pt}%
    \put(0,0){\includegraphics[width=\unitlength]{twist_right_I.eps}}%
    \put(0.46376488,0.07903609){\color[rgb]{0,0,0}\makebox(0,0)[lt]{\lineheight{40.54999924}\smash{\begin{tabular}[t]{l}$\alpha_{n+1}$\end{tabular}}}}%
  \end{picture}%
\endgroup%

%% file: twist_left.eps_tex
\begingroup%
  \makeatletter%
  \providecommand\color[2][]{%
    \errmessage{(Inkscape) Color is used for the text in Inkscape, but the package 'color.sty' is not loaded}%
    \renewcommand\color[2][]{}%
  }%
  \providecommand\transparent[1]{%
    \errmessage{(Inkscape) Transparency is used (non-zero) for the text in Inkscape, but the package 'transparent.sty' is not loaded}%
    \renewcommand\transparent[1]{}%
  }%
  \providecommand\rotatebox[2]{#2}%
  \newcommand*\fsize{\dimexpr\f@size pt\relax}%
  \newcommand*\lineheight[1]{\fontsize{\fsize}{#1\fsize}\selectfont}%
  \ifx\svgwidth\undefined%
    \setlength{\unitlength}{162.63221819bp}%
    \ifx\svgscale\undefined%
      \relax%
    \else%
      \setlength{\unitlength}{\unitlength * \real{\svgscale}}%
    \fi%
  \else%
    \setlength{\unitlength}{\svgwidth}%
  \fi%
  \global\let\svgwidth\undefined%
  \global\let\svgscale\undefined%
  \makeatother%
  \begin{picture}(1,0.42578734)%
    \lineheight{1}%
    \setlength\tabcolsep{0pt}%
    \put(0,0){\includegraphics[width=\unitlength]{twist_left.eps}}%
  \end{picture}%
\endgroup%

%% file: twist_right_II.eps_tex
\begingroup%
  \makeatletter%
  \providecommand\color[2][]{%
    \errmessage{(Inkscape) Color is used for the text in Inkscape, but the package 'color.sty' is not loaded}%
    \renewcommand\color[2][]{}%
  }%
  \providecommand\transparent[1]{%
    \errmessage{(Inkscape) Transparency is used (non-zero) for the text in Inkscape, but the package 'transparent.sty' is not loaded}%
    \renewcommand\transparent[1]{}%
  }%
  \providecommand\rotatebox[2]{#2}%
  \newcommand*\fsize{\dimexpr\f@size pt\relax}%
  \newcommand*\lineheight[1]{\fontsize{\fsize}{#1\fsize}\selectfont}%
  \ifx\svgwidth\undefined%
    \setlength{\unitlength}{162.4427608bp}%
    \ifx\svgscale\undefined%
      \relax%
    \else%
      \setlength{\unitlength}{\unitlength * \real{\svgscale}}%
    \fi%
  \else%
    \setlength{\unitlength}{\svgwidth}%
  \fi%
  \global\let\svgwidth\undefined%
  \global\let\svgscale\undefined%
  \makeatother%
  \begin{picture}(1,0.96339201)%
    \lineheight{1}%
    \setlength\tabcolsep{0pt}%
    \put(0,0){\includegraphics[width=\unitlength]{twist_right_II.eps}}%
    \put(0.42946919,0.07658629){\color[rgb]{0,0,0}\makebox(0,0)[lt]{\lineheight{40.54999924}\smash{\begin{tabular}[t]{l}$\alpha_{n+1}$\end{tabular}}}}%
  \end{picture}%
\endgroup%

%% file: tube_iso_1.eps_tex
\begingroup%
  \makeatletter%
  \providecommand\color[2][]{%
    \errmessage{(Inkscape) Color is used for the text in Inkscape, but the package 'color.sty' is not loaded}%
    \renewcommand\color[2][]{}%
  }%
  \providecommand\transparent[1]{%
    \errmessage{(Inkscape) Transparency is used (non-zero) for the text in Inkscape, but the package 'transparent.sty' is not loaded}%
    \renewcommand\transparent[1]{}%
  }%
  \providecommand\rotatebox[2]{#2}%
  \newcommand*\fsize{\dimexpr\f@size pt\relax}%
  \newcommand*\lineheight[1]{\fontsize{\fsize}{#1\fsize}\selectfont}%
  \ifx\svgwidth\undefined%
    \setlength{\unitlength}{225.66000499bp}%
    \ifx\svgscale\undefined%
      \relax%
    \else%
      \setlength{\unitlength}{\unitlength * \real{\svgscale}}%
    \fi%
  \else%
    \setlength{\unitlength}{\svgwidth}%
  \fi%
  \global\let\svgwidth\undefined%
  \global\let\svgscale\undefined%
  \makeatother%
  \begin{picture}(1,0.63123681)%
    \lineheight{1}%
    \setlength\tabcolsep{0pt}%
    \put(0,0){\includegraphics[width=\unitlength]{tube_iso_1.eps}}%
  \end{picture}%
\endgroup%

%% file: tube_iso_2.eps_tex
\begingroup%
  \makeatletter%
  \providecommand\color[2][]{%
    \errmessage{(Inkscape) Color is used for the text in Inkscape, but the package 'color.sty' is not loaded}%
    \renewcommand\color[2][]{}%
  }%
  \providecommand\transparent[1]{%
    \errmessage{(Inkscape) Transparency is used (non-zero) for the text in Inkscape, but the package 'transparent.sty' is not loaded}%
    \renewcommand\transparent[1]{}%
  }%
  \providecommand\rotatebox[2]{#2}%
  \newcommand*\fsize{\dimexpr\f@size pt\relax}%
  \newcommand*\lineheight[1]{\fontsize{\fsize}{#1\fsize}\selectfont}%
  \ifx\svgwidth\undefined%
    \setlength{\unitlength}{225.14076716bp}%
    \ifx\svgscale\undefined%
      \relax%
    \else%
      \setlength{\unitlength}{\unitlength * \real{\svgscale}}%
    \fi%
  \else%
    \setlength{\unitlength}{\svgwidth}%
  \fi%
  \global\let\svgwidth\undefined%
  \global\let\svgscale\undefined%
  \makeatother%
  \begin{picture}(1,0.67751686)%
    \lineheight{1}%
    \setlength\tabcolsep{0pt}%
    \put(0,0){\includegraphics[width=\unitlength]{tube_iso_2.eps}}%
  \end{picture}%
\endgroup%

%% file: tube_iso_3.eps_tex
\begingroup%
  \makeatletter%
  \providecommand\color[2][]{%
    \errmessage{(Inkscape) Color is used for the text in Inkscape, but the package 'color.sty' is not loaded}%
    \renewcommand\color[2][]{}%
  }%
  \providecommand\transparent[1]{%
    \errmessage{(Inkscape) Transparency is used (non-zero) for the text in Inkscape, but the package 'transparent.sty' is not loaded}%
    \renewcommand\transparent[1]{}%
  }%
  \providecommand\rotatebox[2]{#2}%
  \newcommand*\fsize{\dimexpr\f@size pt\relax}%
  \newcommand*\lineheight[1]{\fontsize{\fsize}{#1\fsize}\selectfont}%
  \ifx\svgwidth\undefined%
    \setlength{\unitlength}{235.6148071bp}%
    \ifx\svgscale\undefined%
      \relax%
    \else%
      \setlength{\unitlength}{\unitlength * \real{\svgscale}}%
    \fi%
  \else%
    \setlength{\unitlength}{\svgwidth}%
  \fi%
  \global\let\svgwidth\undefined%
  \global\let\svgscale\undefined%
  \makeatother%
  \begin{picture}(1,0.48888887)%
    \lineheight{1}%
    \setlength\tabcolsep{0pt}%
    \put(0,0){\includegraphics[width=\unitlength]{tube_iso_3.eps}}%
    \put(0.55544721,0.12666063){\color[rgb]{0,0,0}\makebox(0,0)[lt]{\lineheight{40.54999924}\smash{\begin{tabular}[t]{l}$\alpha_{n+2}$\end{tabular}}}}%
    \put(0.43625208,0.3934305){\color[rgb]{0,0,0}\makebox(0,0)[lt]{\lineheight{40.54999924}\smash{\begin{tabular}[t]{l}$\alpha_{n+1}$\end{tabular}}}}%
  \end{picture}%
\endgroup%

%% file: plus_band_cross_2.eps_tex
\begingroup%
  \makeatletter%
  \providecommand\color[2][]{%
    \errmessage{(Inkscape) Color is used for the text in Inkscape, but the package 'color.sty' is not loaded}%
    \renewcommand\color[2][]{}%
  }%
  \providecommand\transparent[1]{%
    \errmessage{(Inkscape) Transparency is used (non-zero) for the text in Inkscape, but the package 'transparent.sty' is not loaded}%
    \renewcommand\transparent[1]{}%
  }%
  \providecommand\rotatebox[2]{#2}%
  \ifx\svgwidth\undefined%
    \setlength{\unitlength}{199.52443652bp}%
    \ifx\svgscale\undefined%
      \relax%
    \else%
      \setlength{\unitlength}{\unitlength * \real{\svgscale}}%
    \fi%
  \else%
    \setlength{\unitlength}{\svgwidth}%
  \fi%
  \global\let\svgwidth\undefined%
  \global\let\svgscale\undefined%
  \makeatother%
  \begin{picture}(1,0.99123228)%
    \put(0,0){\includegraphics[width=\unitlength]{plus_band_cross_2.eps}}%
  \end{picture}%
\endgroup%

%% file: minus_band_cross_2.eps_tex
\begingroup%
  \makeatletter%
  \providecommand\color[2][]{%
    \errmessage{(Inkscape) Color is used for the text in Inkscape, but the package 'color.sty' is not loaded}%
    \renewcommand\color[2][]{}%
  }%
  \providecommand\transparent[1]{%
    \errmessage{(Inkscape) Transparency is used (non-zero) for the text in Inkscape, but the package 'transparent.sty' is not loaded}%
    \renewcommand\transparent[1]{}%
  }%
  \providecommand\rotatebox[2]{#2}%
  \ifx\svgwidth\undefined%
    \setlength{\unitlength}{199.52442916bp}%
    \ifx\svgscale\undefined%
      \relax%
    \else%
      \setlength{\unitlength}{\unitlength * \real{\svgscale}}%
    \fi%
  \else%
    \setlength{\unitlength}{\svgwidth}%
  \fi%
  \global\let\svgwidth\undefined%
  \global\let\svgscale\undefined%
  \makeatother%
  \begin{picture}(1,0.99123233)%
    \put(0,0){\includegraphics[width=\unitlength]{minus_band_cross_2.eps}}%
  \end{picture}%
\endgroup%

%% file: 3_7_1.eps_tex
\begingroup%
  \makeatletter%
  \providecommand\color[2][]{%
    \errmessage{(Inkscape) Color is used for the text in Inkscape, but the package 'color.sty' is not loaded}%
    \renewcommand\color[2][]{}%
  }%
  \providecommand\transparent[1]{%
    \errmessage{(Inkscape) Transparency is used (non-zero) for the text in Inkscape, but the package 'transparent.sty' is not loaded}%
    \renewcommand\transparent[1]{}%
  }%
  \providecommand\rotatebox[2]{#2}%
  \newcommand*\fsize{\dimexpr\f@size pt\relax}%
  \newcommand*\lineheight[1]{\fontsize{\fsize}{#1\fsize}\selectfont}%
  \ifx\svgwidth\undefined%
    \setlength{\unitlength}{311.49288661bp}%
    \ifx\svgscale\undefined%
      \relax%
    \else%
      \setlength{\unitlength}{\unitlength * \real{\svgscale}}%
    \fi%
  \else%
    \setlength{\unitlength}{\svgwidth}%
  \fi%
  \global\let\svgwidth\undefined%
  \global\let\svgscale\undefined%
  \makeatother%
  \begin{picture}(1,0.9021479)%
    \lineheight{1}%
    \setlength\tabcolsep{0pt}%
    \put(0,0){\includegraphics[width=\unitlength]{3_7_1.eps}}%
  \end{picture}%
\endgroup%

%% file: 3_7_2.eps_tex
\begingroup%
  \makeatletter%
  \providecommand\color[2][]{%
    \errmessage{(Inkscape) Color is used for the text in Inkscape, but the package 'color.sty' is not loaded}%
    \renewcommand\color[2][]{}%
  }%
  \providecommand\transparent[1]{%
    \errmessage{(Inkscape) Transparency is used (non-zero) for the text in Inkscape, but the package 'transparent.sty' is not loaded}%
    \renewcommand\transparent[1]{}%
  }%
  \providecommand\rotatebox[2]{#2}%
  \newcommand*\fsize{\dimexpr\f@size pt\relax}%
  \newcommand*\lineheight[1]{\fontsize{\fsize}{#1\fsize}\selectfont}%
  \ifx\svgwidth\undefined%
    \setlength{\unitlength}{311.15610515bp}%
    \ifx\svgscale\undefined%
      \relax%
    \else%
      \setlength{\unitlength}{\unitlength * \real{\svgscale}}%
    \fi%
  \else%
    \setlength{\unitlength}{\svgwidth}%
  \fi%
  \global\let\svgwidth\undefined%
  \global\let\svgscale\undefined%
  \makeatother%
  \begin{picture}(1,0.90289476)%
    \lineheight{1}%
    \setlength\tabcolsep{0pt}%
    \put(0,0){\includegraphics[width=\unitlength]{3_7_2.eps}}%
  \end{picture}%
\endgroup%

%% file: gl_pairing_final_II.bbl
\providecommand{\bysame}{\leavevmode\hbox to3em{\hrulefill}\thinspace}
\providecommand{\MR}{\relax\ifhmode\unskip\space\fi MR }
\providecommand{\MRhref}[2]{%
  \href{http://www.ams.org/mathscinet-getitem?mr=#1}{#2}
}
\providecommand{\href}[2]{#2}
\begin{thebibliography}{10}

\bibitem{akbulut_kirby}
S.~Akbulut and R.~Kirby, \emph{Branched covers of surfaces in {$4$}-manifolds},
  Math. Ann. \textbf{252} (1979/80), no.~2, 111--131. \MR{593626}

\bibitem{batson}
J.~Batson, \emph{Nonorientable slice genus can be arbitrarily large}, Math.
  Res. Lett. \textbf{21} (2014), no.~3, 423--436. \MR{3272020}

\bibitem{bcg}
H.~U. Boden, M.~Chrisman, and R.~Gaudreau, \emph{Signature and concordance of
  virtual knots}, Indiana Univ. Math. J. \textbf{69} (2020), no.~7, 2395--2459.
  \MR{4195608}

\bibitem{bck}
H.~U. Boden, M.~Chrisman, and H.~Karimi, \emph{The {G}ordon-{L}itherland
  pairing for links in thickened surfaces}, Internat. J. Math. \textbf{33}
  (2022), no.~10-11, Paper No. 2250078, 47. \MR{4514298}

\bibitem{bz}
G.~Burde, H.~Zieschang, and M.~Heusener, \emph{Knots}, extended ed., De Gruyter
  Studies in Mathematics, vol.~5, De Gruyter, Berlin, 2014. \MR{3156509}

\bibitem{cimasoni_turaev}
D.~Cimasoni and V.~Turaev, \emph{A generalization of several classical
  invariants of links}, Osaka J. Math. \textbf{44} (2007), no.~3, 531--561.
  \MR{2360939}

\bibitem{clark_78}
B.~E. Clark, \emph{Crosscaps and knots}, Internat. J. Math. Math. Sci.
  \textbf{1} (1978), no.~1, 113--123. \MR{478131}

\bibitem{conway_sloane}
J.~H. Conway and N.~J.~A. Sloane, \emph{Sphere packings, lattices and groups},
  third ed., Grundlehren der Mathematischen Wissenschaften [Fundamental
  Principles of Mathematical Sciences], vol. 290, ch.~``On the classification
  of integral quadratic forms'', pp.~lxxiv+703, Springer-Verlag, New York,
  1999, With additional contributions by E. Bannai, R. E. Borcherds, J. Leech,
  S. P. Norton, A. M. Odlyzko, R. A. Parker, L. Queen and B. B. Venkov.
  \MR{1662447}

\bibitem{dasbach_lowrance}
O.~T. Dasbach and A.~M. Lowrance, \emph{Turaev genus, knot signature, and the
  knot homology concordance invariants}, Proc. Amer. Math. Soc. \textbf{139}
  (2011), no.~7, 2631--2645. \MR{2784832}

\bibitem{gay_kirby}
D.~Gay and R.~Kirby, \emph{Trisecting 4-manifolds}, Geom. Topol. \textbf{20}
  (2016), no.~6, 3097--3132. \MR{3590351}

\bibitem{gilmer_livingston}
P.~M. Gilmer and C.~Livingston, \emph{The nonorientable 4-genus of knots}, J.
  Lond. Math. Soc. (2) \textbf{84} (2011), no.~3, 559--577. \MR{2855790}

\bibitem{gompf_stipsicz}
R.~E. Gompf and A.~I. Stipsicz, \emph{{$4$}-manifolds and {K}irby calculus},
  Graduate Studies in Mathematics, vol.~20, American Mathematical Society,
  Providence, RI, 1999. \MR{1707327}

\bibitem{gordon_litherland}
C.~McA. Gordon and R.~A. Litherland, \emph{On the signature of a link}, Invent.
  Math. \textbf{47} (1978), no.~1, 53--69. \MR{500905}

\bibitem{greene}
J.~E. Greene, \emph{Alternating links and definite surfaces}, Duke Math. J.
  \textbf{166} (2017), no.~11, 2133--2151, With an appendix by Andr\'{a}s
  Juh\'{a}sz and Marc Lackenby. \MR{3694566}

\bibitem{hsvc}
J.~Hoste, P.~D. Shanahan, and C.~A. Van~Cott, \emph{On the nonorientable
  4-genus of double twist knots}, J. Knot Theory Ramifications \textbf{32}
  (2023), no.~5, Paper No. 2350036, 30. \MR{4598114}

\bibitem{ILL}
Y.~H. Im, K.~Lee, and S.~Y. Lee, \emph{Signature, nullity and determinant of
  checkerboard colorable virtual links}, J. Knot Theory Ramifications
  \textbf{19} (2010), no.~8, 1093--1114. \MR{2718629}

\bibitem{jabuka_vancott}
S.~Jabuka and C.~A. Van~Cott, \emph{Comparing nonorientable three genus and
  nonorientable four genus of torus knots}, J. Knot Theory Ramifications
  \textbf{29} (2020), no.~3, 2050013, 15. \MR{4101607}

\bibitem{JMMZ}
J.~Joseph, J.~Meier, M.~Miller, and A.~Zupan, \emph{Bridge trisections and
  classical knotted surface theory}, Pacific J. Math. \textbf{319} (2022),
  no.~2, 343--369. \MR{4482720}

\bibitem{kalfagianni_lee}
E.~Kalfagianni and C.~R.~S. Lee, \emph{Crosscap numbers and the {J}ones
  polynomial}, Adv. Math. \textbf{286} (2016), 308--337. \MR{3415687}

\bibitem{kanenobu}
T.~Kanenobu, \emph{Sharp-unknotting number of a torus knot}, Kyungpook Math. J.
  \textbf{49} (2009), no.~3, 583--594. \MR{2601857}

\bibitem{kauffman_87}
L.~H. Kauffman, \emph{State models and the {J}ones polynomial}, Topology
  \textbf{26} (1987), no.~3, 395--407. \MR{899057}

\bibitem{kauffman_taylor}
L.~H. Kauffman and L.~R. Taylor, \emph{Signature of links}, Trans. Amer. Math.
  Soc. \textbf{216} (1976), 351--365. \MR{388373}

\bibitem{khovanov}
M.~Khovanov, \emph{A categorification of the {J}ones polynomial}, Duke Math. J.
  \textbf{101} (2000), no.~3, 359--426. \MR{1740682}

\bibitem{lee}
E.~S. Lee, \emph{An endomorphism of the {K}hovanov invariant}, Adv. Math.
  \textbf{197} (2005), no.~2, 554--586. \MR{2173845}

\bibitem{alexander_machon}
T.~Machon and G.~P. Alexander, \emph{Global defect topology in nematic liquid
  crystals}, Proc. R. Soc. A \textbf{472} (2016), 20160265.

\bibitem{massey}
W.~S. Massey, \emph{Proof of a conjecture of {W}hitney}, Pacific J. Math.
  \textbf{31} (1969), 143--156. \MR{250331}

\bibitem{meier_zupan_2}
J.~Meier and A.~Zupan, \emph{Bridge trisections of knotted surfaces in
  4-manifolds}, Proc. Natl. Acad. Sci. USA \textbf{115} (2018), no.~43,
  10880--10886. \MR{3871791}

\bibitem{murakami}
H.~Murakami, \emph{Some metrics on classical knots}, Math. Ann. \textbf{270}
  (1985), no.~1, 35--45. \MR{769605}

\bibitem{murakami_sakai}
H.~Murakami and S.~Sakai, \emph{Sharp-unknotting number and the {A}lexander
  module}, Topology Appl. \textbf{52} (1993), no.~2, 169--179. \MR{1241192}

\bibitem{murakami_yasuhara}
H.~Murakami and A.~Yasuhara, \emph{Crosscap number of a knot}, Pacific J. Math.
  \textbf{171} (1995), no.~1, 261--273. \MR{1362987}

\bibitem{murasugi_87}
K.~Murasugi, \emph{Jones polynomials and classical conjectures in knot theory},
  Topology \textbf{26} (1987), no.~2, 187--194. \MR{895570}

\bibitem{nicholson}
V.~A. Nicholson, \emph{Twisted surfaces}, J. Math. Chem. \textbf{40} (2006),
  no.~2, 105--117. \MR{2256682}

\bibitem{os}
P.~Ozsv\'{a}th and Z.~Szab\'{o}, \emph{Holomorphic disks and knot invariants},
  Adv. Math. \textbf{186} (2004), no.~1, 58--116. \MR{2065507}

\bibitem{rasmussen_thesis}
J.~A. Rasmussen, \emph{Floer homology and knot complements}, ProQuest LLC, Ann
  Arbor, MI, 2003, Thesis (Ph.D.)--Harvard University. \MR{2704683}

\bibitem{shumakovitch_1}
A.~N. Shumakovitch, \emph{Torsion of {K}hovanov homology}, Fund. Math.
  \textbf{225} (2014), no.~1, 343--364. \MR{3205577}

\bibitem{shumakovitch_2}
\bysame, \emph{Torsion in {K}hovanov homology of homologically thin knots}, J.
  Knot Theory Ramifications \textbf{30} (2021), no.~14, Paper No. 2141015, 17.
  \MR{4407084}

\bibitem{thistlethwaite_87}
M.~B. Thistlethwaite, \emph{A spanning tree expansion of the {J}ones
  polynomial}, Topology \textbf{26} (1987), no.~3, 297--309. \MR{899051}

\bibitem{turaev_87}
V.~G. Turaev, \emph{A simple proof of the {M}urasugi and {K}auffman theorems on
  alternating links}, Enseign. Math. (2) \textbf{33} (1987), no.~3-4, 203--225.
  \MR{925987}

\bibitem{viro}
O.~Ja. Viro, \emph{Positioning in codimension {$2$}, and the boundary}, Uspehi
  Mat. Nauk \textbf{30} (1975), no.~1(181), 231--232. \MR{420641}

\bibitem{yasuhara_96}
A.~Yasuhara, \emph{Connecting lemmas and representing homology classes of
  simply connected {$4$}-manifolds}, Tokyo J. Math. \textbf{19} (1996), no.~1,
  245--261. \MR{1391941}

\bibitem{yasuhara}
\bysame, \emph{An elementary proof for that all unoriented spanning surfaces of
  a link are related by attaching/deleting tubes and {M}\"{o}bius bands}, J.
  Knot Theory Ramifications \textbf{23} (2014), no.~1, 1450004, 5. \MR{3190128}

\end{thebibliography}
